\documentclass[12pt]{amsart}

\usepackage{amsmath,amssymb,amsfonts,amsthm,latexsym,graphicx,multirow}
\usepackage{hyperref}
\usepackage[all]{xy}

\newcommand\A{\mathrm{A}}   \newcommand\Aut{\mathrm{Aut}}
  \newcommand\bighat{\text{\scalebox{1.5}[1]{$\hat{~}$}}}
\newcommand\Cen{\mathbf{C}}  \newcommand\calC{\mathcal{C}} \newcommand\calO{\mathcal{O}}   \newcommand\Co{\mathrm{Co}} 
\newcommand\D{\mathrm{D}} \newcommand\di{\,\big|\,}

\newcommand\G{\mathrm{G}}  \newcommand\GaSp{\mathrm{\Gamma Sp}} \newcommand\GaO{\mathrm{\Gamma O}}  \newcommand\GL{\mathrm{GL}} \newcommand\GO{\mathrm{O}} \newcommand\GU{\mathrm{GU}}

\newcommand\M{\mathrm{M}} \newcommand\magma{{\sc Magma} }
\newcommand\N{\mathrm{N}}
\newcommand\Nor{\mathbf{N}}
\newcommand\Out{\mathrm{Out}}
\newcommand\Pa{\mathrm{P}}     \newcommand\ppd{\mathrm{ppd}} \newcommand\POm{\mathrm{P\Omega}} \newcommand\PSL{\mathrm{PSL}}    \newcommand\PSp{\mathrm{PSp}} \newcommand\PSU{\mathrm{PSU}}

 \newcommand\Rad{\mathrm{Rad}}  
\newcommand\SL{\mathrm{SL}}  \newcommand\Soc{\mathrm{Soc}} \newcommand\Sp{\mathrm{Sp}} \newcommand\Stab{\mathrm{Stab}}   \newcommand\Sy{\mathrm{S}}  \newcommand\Sz{\mathrm{Sz}}

\newtheorem{theorem}{Theorem}[section]
\newtheorem{lemma}[theorem]{Lemma}

\newtheorem{corollary}[theorem]{Corollary}

\theoremstyle{definition}

\newtheorem{example}[theorem]{Example}
\newtheorem*{remark}{Remark}
\newtheorem{hypothesis}[theorem]{Hypothesis}

\begin{document}

\title[Factorizations of almost simple groups]{Factorizations of almost simple groups with a factor having many nonsolvable composition factors}

\author[Li]{Cai Heng Li}
\address{(Li) Department of Mathematics, Southern University of Science and Technology\\Shenzhen, Guangdong 518055\\P. R. China}
\email{lich@sustech.edu.cn}

\author[Xia]{Binzhou Xia}
\address{(Xia) School of Mathematics and Statistics\\The University of Melbourne\\Parkville, VIC 3010\\ Australia}
\email{binzhoux@unimelb.edu.au}

\maketitle

\begin{abstract}
This paper classifies the factorizations of almost simple groups with a factor having at least two nonsolvable composition factors. This together with a previous classification result of the authors reduces the factorization problem of almost simple groups to the case where both factors have a unique nonsolvable composition factor.


\end{abstract}

\section{Introduction}

For a group $G$, an expression $G=HK$ with subgroups $H$ and $K$ of $G$ is called a \emph{factorization} of $G$, where $H$ and $K$ are called  \emph{factors}. Factorizations of groups play an important role not only in group theory but also in other branches of mathematics such as Galois theory~\cite{FGS1993,GS1995} and graph theory~\cite{LX,LPS2010}.

A group $G$ is \emph{almost simple} with \emph{socle} $L$ if $L\leqslant G\leqslant\Aut(L)$ for some finite nonabelian simple group $L$. To study factorizations of finite groups, it is crucial to study those of almost simple groups. In 1987, Hering, Liebeck and Saxl~\cite{HLS1987} classified factorizations of exceptional groups of Lie type. A factorization $G=HK$ is called a \emph{maximal factorization} of $G$ if both $H$ and $K$ are maximal subgroups of $G$. In 1990, Liebeck, Praeger and Saxl published the landmark work~\cite{LPS1990} classifying maximal factorizations of almost simple groups. Furthermore, when the socle is an alternating group, all the factorizations of such a group were determined in~\cite[Theorem~D]{LPS1990}. Based on the maximal factorizations in~\cite[Theorem~C]{LPS1990}, Giudici~\cite{Giudici2006} in 2006 determined the factorizations of sporadic groups. In 2010, all factorizations of almost simple groups for which one factor is maximal and the intersection of the two factors is trivial have been determined~\cite{LPS2010} by Liebeck, Praeger and Saxl.

In~\cite{LX}, factorizations of almost simple groups with a solvable factor are classified. As a counter part, we classify in this paper (Theorem~\ref{thm1}) factorizations of almost simple groups with a factor having at least two nonsolvable composition factors (counted with multiplicity). These two classification results then reduce the factorization problem of almost simple groups to the case where both factors have a unique nonsolvable composition factor (Corollary~\ref{cor1}).

\begin{theorem}\label{thm1}
Let $G$ be an almost simple group with socle $L$. Suppose that $G$ has a factorization $G=HK$, where $H$ has at least two nonsolvable composition factors and $K$ is core-free in $G$. Then either~\emph{(a)} or~\emph{(b)} below holds.
\begin{itemize}
\item[(a)] $\A_n\leqslant G\leqslant\Sy_n$ with $n\geqslant10$, and one of the following holds:
\begin{itemize}
\item[(a.1)] $H$ is a transitive permutation group of degree $n$, and $\A_{n-1}\leqslant K\leqslant\Sy_{n-1}$;
\item[(a.2)] $n=10$, $H$ is a transitive permutation group of degree $10$ such that $(\A_5\times\A_5).2\leqslant H\leqslant\Sy_5\wr\Sy_2$, and $K=\SL_2(8)$ or $\SL_2(8).3$;
\item[(a.3)] $n=12$, $\A_7\times\A_5\leqslant H\leqslant\Sy_7\times\Sy_5$, and $K=\M_{12}$;
\item[(a.4)] $n=24$, $\A_{19}\times\A_5\leqslant H\leqslant\Sy_{19}\times\Sy_5$, and $K=\M_{24}$.
\end{itemize}
\item[(b)] $L$ is a symplectic group or an orthogonal group of plus type, and the triple $(L,H\cap L,K\cap L)$ lies in \emph{Table~\ref{tab1}} or \emph{Table~\ref{tab5}}, respectively, where $Q\leqslant2$.
\end{itemize}
\end{theorem}

{\small
\begin{table}[htbp]
\caption{Factorizations of symplectic groups in Theorem~\ref{thm1}(b)}\label{tab1}
\centering
\begin{tabular}{|l|l|l|l|l|}
\hline
row~& $L$ & $H\cap L$ & $K\cap L$ & Ex.\\
\hline
\multirow{2}*{1} & \multirow{2}*{$\Sp_{4\ell}(2^f)$, $f\ell\geqslant2$} & $(\Sp_{2a}(2^{fb})\times\Sp_{2a}(2^{fb})).R.2$, & \multirow{2}*{$\Omega_{4\ell}^-(2^f).Q$} & \multirow{2}*{\ref{exa1}}\\
 &  & $ab=\ell$ and $R\leqslant b\times b$ &  & \\
\hline
\multirow{2}*{2} & \multirow{2}*{$\Sp_{12\ell}(2^f)$} & $(\G_2(2^{f\ell})\times\G_2(2^{f\ell})).R.2$, & \multirow{2}*{$\Omega_{12\ell}^-(2^f).Q$} &
\multirow{2}*{\ref{exa9}}\\
 &  & $R\leqslant\ell\times\ell$ &  & \\
\hline
3 & $\Sp_{4\ell}(4)$, $\ell\geqslant2$ & $\Sp_2(4)\times\Sp_{4\ell-2}(4)$ & $\Sp_{2\ell}(16).Q$ & max\\
\hline
4 & $\Sp_{4\ell}(4)$, $\ell\geqslant2$ & $(\Sp_2(4)\times\Sp_{2\ell}(4)).P$, $P\leqslant2$ & $\Omega_{4\ell}^-(4).Q$ & \ref{exa2}\\
\hline
5 & $\Sp_{8\ell}(2)$ & $(\Sp_2(4)\times\Sp_{2\ell}(4)).P$, $P\leqslant[8]$ & $\Omega_{8\ell}^-(2).Q$ & \ref{exa3}\\
\hline
6 & $\Sp_4(2^f)$, $f\geqslant3$ odd & $(\Sp_2(2^f)\times\Sp_2(2^f)).P$, $P\leqslant2$ & $\Sz(2^f)$ & max\\
\hline
\multirow{2}*{7} & \multirow{2}*{$\Sp_6(2^f)$, $f\geqslant2$} & $U\times\Sp_4(2^f)$, & \multirow{2}*{$\G_2(2^f)$} & \multirow{2}*{\ref{exa4}}\\
 &  & $U$ nonsolvable and $U\leqslant\Sp_2(2^f)$ &  & \\
\hline
8 & $\Sp_{12}(4)$ & $\Sp_2(4)\times\Sp_{10}(4)$ & $\G_2(16).Q$ & \ref{exa12}\\
\hline
9 & $\Sp_{12}(4)$ & $(\Sp_2(4)\times\G_2(4)).P$, $P\leqslant2$ & $\Omega_{12}^-(4).Q$ & \ref{exa13}\\
\hline
10 & $\Sp_{24}(2)$ & $(\Sp_2(4)\times\G_2(4)).P$, $P\leqslant[8]$ & $\Omega_{24}^-(2).Q$ & \ref{exa15}\\
\hline
\end{tabular}
\end{table}

\begin{table}[htbp]
\caption{Factorizations of orthogonal groups in Theorem~\ref{thm1}(b)}\label{tab5}
\centering
\begin{tabular}{|l|l|l|l|l|}
\hline
row~& $L$ & $H\cap L$ & $K\cap L$ & Ex.\\
\hline
\multirow{2}*{1} & \multirow{2}*{$\POm_{4\ell}^+(q)$, $\ell\geqslant2$, $q\geqslant4$} & $(U\times\PSp_{2\ell}(q)).P$, $U$ nonsolvable, & \multirow{2}*{$\Omega_{4\ell-1}(q)$} & \multirow{2}*{\ref{exa5}}\\
 &  & $U\leqslant\PSp_2(q)$ and $P\leqslant\gcd(2,\ell,q-1)$ & & \\
\hline
2 & $\Omega_{8\ell}^+(2)$ & $(\Sp_2(4)\times\Sp_{2\ell}(4)).P$, $P\leqslant2^2$ & $\Sp_{8\ell-2}(2)$ & \ref{exa10}\\
\hline
\multirow{2}*{3} & \multirow{2}*{$\Omega_{8\ell}^+(4)$} & $(\Sp_2(4^c)\times\Sp_{2\ell}(16)).P$, & \multirow{2}*{$\Sp_{8\ell-2}(4)$} & \multirow{2}*{\ref{exa6}}\\
 &  & $c=1$ or $2$ and $P\leqslant2^2$ &  & \\
\hline
\multirow{2}*{4} & \multirow{2}*{$\Omega_{12}^+(2^f)$, $f\geqslant2$} & $U\times\G_2(2^f)$, & \multirow{2}*{$\Sp_{10}(2^f)$} & \multirow{2}*{\ref{exa11}}\\
 &  & $U$ nonsolvable and $U\leqslant\Sp_2(2^f)$ & & \\
\hline
5 & $\Omega_{24}^+(2)$ & $(\Sp_2(4)\times\G_2(4)).P$, $P\leqslant2^2$ & $\Sp_{22}(4)$ & \ref{exa14}\\
\hline
\multirow{2}*{6} & \multirow{2}*{$\Omega_{24}^+(4)$} & $(\Sp_2(4^c)\times\G_2(16)).P$, & \multirow{2}*{$\Sp_{22}(4)$} & \multirow{2}*{\ref{exa16}}\\
 &  & $c=1$ or $2$ and $P\leqslant2^2$ &  &  \\
\hline
\end{tabular}
\end{table}
}

Here are some remarks on Theorem~\ref{thm1}.
\begin{itemize}
\item[(i)] For each factorization $G=HK$ described in Theorem~\ref{thm1}, $K\cap L$ is almost simple, and either $H$ has exactly two nonsolvable composition factors or~(a.1) occurs.
\item[(ii)] It is easy to see that each case in part~(a) of Theorem~\ref{thm1} gives rise to a factorization $G=HK$. Also, for each row of Table~\ref{tab1} and Table~\ref{tab5} there is an example for such factorization $G=HK$, which can be found in the last column: the label $3.k$ indicates that such an example is given in Example~$3.k$ while the label max indicates that there exists a maximal factorization (see Lemma~\ref{lem1}).
\end{itemize}

\begin{corollary}\label{cor1}
Let $G$ be an almost simple group. Suppose $G=HK$ with core-free subgroups $H$ and $K$ of $G$. Then interchanging $H$ and $K$ if necessary, one of the following holds.
\begin{itemize}
\item[(a)] $H$ is solvable, and $(G,H,K)$ is described in~\cite[Theorem~1.1]{LX}.
\item[(b)] $H$ has at least two nonsolvable composition factors, and $(G,H,K)$ is described in \emph{Theorem~\ref{thm1}}.
\item[(c)] Both $H$ and $K$ have a unique nonsolvable composition factor.
\end{itemize}
\end{corollary}

The proof of Theorem~\ref{thm1} is sketched in Section~4 and completed in Sections~5 and~6. In Section~3, we construct factorizations as described in Table~\ref{tab1} and Table~\ref{tab5}.

\vskip0.1in
\noindent\textsc{Acknowledgements.}
This work was supported by NSFC grants 11231008 and 11771200. The second author is very grateful to Southern University of Science and Technology for the financial support of his visit to the first author. The authors would like to thank the anonymous referee for very helpful comments.

\section{Preliminaries}

All the groups in this paper are assumed to be finite. Our group-theoretic notation mostly follows \cite{LPS1990}.

For a positive integer $m$ and prime number $p$, denote by $m_p$ the largest $p$-power that divides $m$. Given positive integers $a$ and $n$, a prime number $r$ is called a \emph{primitive prime divisor} of the pair $(a,n)$ if $r$ divides $a^n-1$ but does not divide $a^i-1$ for any positive integer $i<n$. In other words, a primitive prime divisor of $(a,n)$ is a prime number $r$ such that $a$ has order $n$ in $\mathbb{F}_r^\times$. In particular, we have the observation in the next lemma, which will be used repeatedly (and sometimes implicitly) in this paper.

\begin{lemma}\label{lem7}
If $r$ is a primitive prime divisor of $(a,m)$, then $m\di r-1$ and so $r>m$.
\end{lemma}

By an elegant theorem of Zsigmondy (see for example~\cite[Theorem IX.8.3]{Blackburn1982}), $(a,n)$ has a primitive prime divisor whenever $a\geqslant2$ and $n\geqslant3$ with $(a,n)\neq(2,6)$. For positive integers $a\geqslant2$ and $n\geqslant3$, denote the set of primitive prime divisors of $(a,n)$ by $\ppd(a,n)$ if $(a,n)\neq(2,6)$, and set $\ppd(2,6)=\{7\}$.

For a group $X$ and a subgroup $Y$, denote the set of right cosets of $Y$ in $X$ by $[X:Y]$. There are several equivalent conditions for a group factorization as in Lemma~\ref{lem6} below.

\begin{lemma}\label{lem6}
Let $H$ and $K$ be subgroups of $G$. Then the following are equivalent.
\begin{itemize}
\item[(a)] $G=HK$.
\item[(b)] $G=H^xK^y$ for any $x,y\in G$.
\item[(c)] $|H\cap K||G|=|H||K|$.
\item[(d)] $|G|\leqslant|H||K|/|H\cap K|$.
\item[(e)] $H$ acts transitively on $[G:K]$ by right multiplication.
\item[(f)] $K$ acts transitively on $[G:H]$ by right multiplication.
\end{itemize}
\end{lemma}

\begin{remark}\label{rmk2}
Lemma~\ref{lem6} is easy to prove but plays a fundamental role in the study of group factorizations. For example,
\begin{itemize}
\item[(i)] due to part~(b) we will consider conjugacy classes of subgroups when studying factorizations of a group;
\item[(ii)] given a group $G$ and its subgroups $H$ and $K$, in order to inspect whether $G=HK$ we only need to compute the orders of $G$, $H$, $K$ and $H\cap K$ by part~(c) or~(d), which enables us to search factorizations of a group efficiently in \magma~\cite{magma}.
\end{itemize}
\end{remark}

Here are more observations on group factorizations.

\begin{lemma}\label{lem8}
Let $H$ and $K$ be subgroups of $G$, and $L$ be a normal subgroup of $G$. If $G=HK$, then we have the following divisibilities.
\begin{itemize}
\item[(a)] $|G|$ divides $|H||K|$.
\item[(b)] $|G|$ divides $|H\cap L||K||G/L|$.
\item[(c)] $|L|$ divides $|H\cap L||K|$.
\item[(d)] $|L|$ divides $|H\cap L||K\cap L||G/L|$.
\end{itemize}
\end{lemma}

\begin{lemma}\label{lem9}
Let $H$, $K$ and $M$ be subgroups of $G$. If $G=HK$ and $H\leqslant M$, then $M=H(K\cap M)$.
\end{lemma}

\begin{lemma}\label{lem16}
Let $H$, $K$ and $M$ be subgroups of $G$. Suppose $H\leqslant M$ and $M=H(K\cap M)$. Then $G=HK$ if and only if $G=MK$.
\end{lemma}

To construct more factorizations from existing ones, we introduce some diagrams for convenience: for a group $G$ and subgroups $H$, $K$ and $M$ of $G$, the diagram
{\large
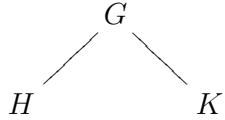
$$
\scalebox{0.8}[0.8]{
\xymatrix{
&G\ar@{-}[dl]\ar@{-}[dr]&\\
H&&K
}
}
$$
}
means $G=HK$, while the diagram
{\large
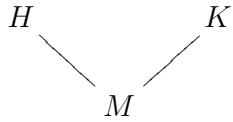
$$
\scalebox{0.8}[0.8]{
\xymatrix{
H\ar@{-}[dr]&&K\ar@{-}[dl]\\
&M&
}
}
$$
}
means $H\cap K=M$. Then the diagram
{\large
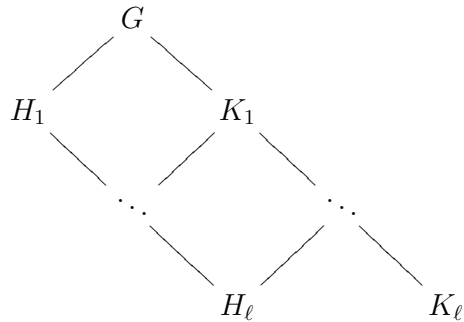
$$
\scalebox{0.8}[0.8]{
\xymatrix{
&&G\ar@{-}[dl]\ar@{}[d]\ar@{-}[dr]&&&&\\
&H_1\ar@{-}[dr]&&K_1\ar@{-}[dl]\ar@{-}[dr]&&\\
&&\ddots\ar@{-}[dr]&&\ddots\ar@{-}[dl]\ar@{-}[dr]\\
&&&H_\ell&&K_\ell&
}
}
$$
}
implies $G=H_1K_1=H_1H_2K_2=\dots=H_1H_2\cdots H_\ell K_\ell=H_1K_\ell$.

\section{Examples}

Let $S$ be a classical group defined on a vector space $V$. If $Z$ is the center of $S$ such that $S/Z$ is simple group and $X$ is a subgroup of $\GL(V)$ containing $S$ as a normal subgroup, then for a subgroup $Y$ of $X$, denote by $\bighat Y$ the subgroup $(Y\cap S)Z/Z$ of $S/Z$. If $S$ is a linear group, define $\Pa_k[S]$ to be the stabilizer of a $k$-space in $S$. Next assume that $S$ is a symplectic, unitary or orthogonal group.
\begin{itemize}
\item[(i)] If $S$ is transitive on the set of totally singular $k$-spaces, then let $\Pa_k[S]$ be the stabilizer of a totally singular $k$-space in $S$.
\item[(ii)] If $S$ is not transitive on the set of totally singular $k$-spaces, then let $\Pa_{k-1}[S]$ and $\Pa_k[S]$ be the stabilizers of totally singular $k$-spaces in the two different orbits of $S$.
\item[(iii)] If $S$ is not transitive on the set of totally singular $k$-spaces, then let $\Pa_{k-1,k}[S]$ be the stabilizer of a totally singular $(k-1)$-space in $S$.
\item[(iv)] If $S$ is not transitive on the set of totally singular $k$-spaces, then let $\Pa_{1,k-1}[S]$ be the intersection $\Pa_1[S]\cap\Pa_{k-1}[S]$, where the $1$-space stabilized by $\Pa_1[S]$ lies in the $k$-space stabilized by $\Pa_{k-1}[S]$.
\item[(v)] If $S$ is not transitive on the set of totally singular $k$-spaces, then let $\Pa_{1,k}[S]$ be the intersection $\Pa_1[S]\cap\Pa_k[S]$, where the $1$-space stabilized by $\Pa_1[S]$ lies in the $k$-space stabilized by $\Pa_k[S]$.
\end{itemize}
For a non-degenerate $k$-space $W$,
\begin{itemize}
\item[(i)] denote the stabilizer of $W$ in $S$ by $\N_k[S]$ if either $S$ is symplectic or unitary, or $S$ is orthogonal of even dimension with $k$ odd;
\item[(ii)] denote the stabilizer of $W$ in $S$ by $\N_k^\varepsilon[S]$ for $\varepsilon=\pm$ if $S$ is orthogonal and $W$ has type $\varepsilon$;
\item[(iii)] denote the stabilizer of $W$ in $S$ by $\N_k^\varepsilon[S]$ for $\varepsilon=\pm$ if $S$ is orthogonal of odd dimension and $W^\bot$ has type~$\varepsilon$.
\end{itemize}
For the above defined groups $\Pa_k[S]$, $\Pa_{i,j}[S]$, $\N_k[S]$, $\N_k^-[S]$ and $\N_k^+[S]$, we will simply write $\Pa_k$, $\Pa_{i,j}$, $\N_k$, $\N_k^-$ and $\N_k^+$, respectively, if the classical group $S$ is clear from the context.

Let $G$ be an almost simple group with socle classical simple, and assume that no element of $G$ induces a triality automorphism if $\Soc(G)=\POm_8^+(q)$. In~\cite{Aschbacher1984}, Aschbacher defined eight families $\calC_1$--$\calC_8$ of subgroups of $G$. These groups are now called \emph{geometric subgroups} of $G$, and described in more detail by Kleidman and Liebeck~\cite{KL1990}.

The following lemma can be read off from Tables~1--4 of~\cite{LPS1990}.

\begin{lemma}\label{lem1}
Let $G$ be an almost simple group with socle $L$ classical. If $G=AB$ is a maximal factorization of $G$ such that $A$ has at least two nonsolvable composition factors and $B$ is core-free in $G$, then $A$ has exactly two nonsolvable composition factors and $(L,A\cap L,B\cap L)=(L,X,Y)$ as in \emph{Table~\ref{tab2}}. Conversely, for each triple $(L,X,Y)$ in \emph{Table~\ref{tab2}}, there exists a maximal factorization $G=HK$ satisfying $\Soc(G)=L$, $H\cap L=X$ and $K\cap L=Y$.
\end{lemma}

\begin{table}[htbp]
\caption{Maximal factorizations of classical almost simple groups with a factor having at least two nonsolvable composition factors}\label{tab2}
\centering
\begin{tabular}{|l|l|l|l|}
\hline
row & $L$ & $X$ & $Y$\\
\hline
1 & $\Sp_{4\ell}(2^f)$, $f\ell\geqslant2$ & $(\Sp_{2\ell}(2^f)\times\Sp_{2\ell}(2^f)).2$ & $\GO_{4\ell}^-(2^f)$\\
2 & $\Sp_{4\ell}(4)$, $\ell\geqslant2$ & $\Sp_2(4)\times\Sp_{4\ell-2}(4)$ & $\Sp_{2\ell}(16).2$\\
3 & $\Sp_4(2^f)$, $f\geqslant3$ odd & $(\Sp_2(2^f)\times\Sp_2(2^f)).2$ & $\Sz(2^f)$\\
4 & $\Sp_6(2^f)$, $f\geqslant2$ & $\Sp_2(2^f)\times\Sp_4(2^f)$ & $\G_2(2^f)$\\
\hline
5 & $\POm_{4\ell}^+(q)$, $\ell\geqslant3$, $q\geqslant4$ & $(\PSp_2(q)\times\PSp_{2\ell}(q)).\gcd(2,\ell,q-1)$ & $\Omega_{4\ell-1}(q)$\\
6 & $\POm_8^+(q)$, $q\geqslant5$ odd & $(\PSp_2(q)\times\PSp_4(q)).2$ & $\Omega_7(q)$\\
7 & $\Omega_8^+(2)$ & $(\SL_2(4)\times\SL_2(4)).2^2$ & $\Sp_6(2)$\\
8 & $\Omega_8^+(4)$ & $(\SL_2(16)\times\SL_2(16)).2^2$ & $\Sp_6(4)$\\
\hline
\end{tabular}
\end{table}

\begin{example}\label{exa1}
Let $G=\Sp_{4\ell}(2^f)$ with $f\ell\geqslant2$, and $K=\GO_{4\ell}^-(2^f)<G$. Then for all positive integers $a$ and $b$ such that $ab=\ell$, there exists a subgroup $H=\Sp_{2a}(2^{fb})\wr\Sy_2=((\Sp_{2a}(2^{fb})\times\Sp_{2a}(2^{fb})){:}2$ of $G$ such that $G=HK$, as in row~1 of Table~\ref{tab1}.
\end{example}

In fact, let $M=\Sp_{4a}(2^{fb}){:}b$ be a $\calC_3$-subgroup of $G$. Then $M\cap K=\GO_{4a}^-(2^{fb}).b$, and so $|M\cap K|||G|=|M||K|$. It follows that $G=MK$ by Lemma~\ref{lem6}, and then we have the diagram
{\large
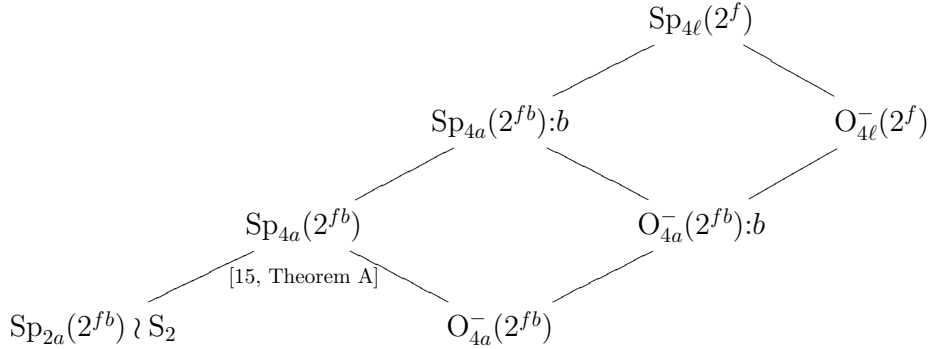
$$
\scalebox{0.8}[0.8]{
\xymatrix{
&&&\Sp_{4\ell}(2^f)\ar@{-}[dl]\ar@{}[d]\ar@{-}[dr]&\\
&&\Sp_{4a}(2^{fb}){:}b\ar@{-}[dl]\ar@{-}[dr]&&\GO_{4\ell}^-(2^f)\ar@{-}[dl]\\
&\Sp_{4a}(2^{fb})\ar@{-}[dl]\ar@{}[d]|{\text{\cite[Theorem~A]{LPS1990}}}\ar@{-}[dr]&&\GO_{4a}^-(2^{fb}){:}b\ar@{-}[dl]\\
\Sp_{2a}(2^{fb})\wr\Sy_2&&\GO_{4a}^-(2^{fb})&
}
}
$$
}
This gives a subgroup $H=\Sp_{2a}(2^{fb})\wr\Sy_2=((\Sp_{2a}(2^{fb})\times\Sp_{2a}(2^{fb})){:}2$ of $G$ such that $G=HK$, proving Example~\ref{exa1}.

\begin{example}\label{exa9}
Let $G=\Sp_{12\ell}(2^f)$ and $K=\GO_{12\ell}^-(2^f)<G$. Then there exists a subgroup $H=G_2(2^{f\ell})\wr\Sy_2=(\G_2(2^f)\times\G_2(2^f)){:}2$ of $G$ such that $G=HK$, as in row~2 of Table~\ref{tab1}.
\end{example}

To show that Example~\ref{exa9} holds, we first note $\Sp_6(2^{f\ell})=\G_2(2^{f\ell})\GO_6^\varepsilon(2^{f\ell})$ for $\varepsilon=\pm$ by~\cite[Theorem~A]{LPS1990}. This implies that
$$
\Sp_6(2^{f\ell})\times\Sp_6(2^{f\ell})=(\G_2(2^{f\ell})\times\G_2(2^{f\ell})(\GO_6^-(2^{f\ell})\times\GO_6^+(2^{f\ell}))
$$
and so $\Sp_6(2^{f\ell})\wr\Sy_2=(\G_2(2^{f\ell})\wr\Sy_2)(\GO_6^-(2^{f\ell})\times\GO_6^+(2^{f\ell}))$. Then the diagram
{\large
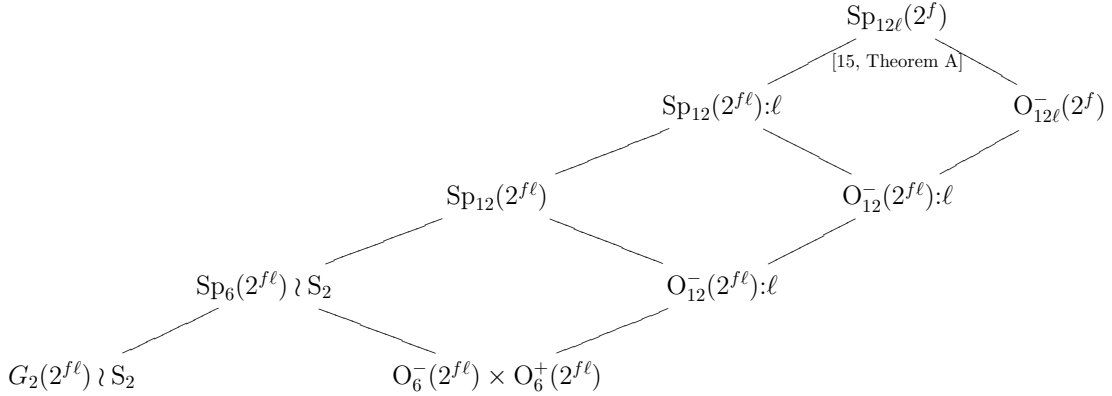
$$
\scalebox{0.7}[0.7]{
\xymatrix{
&&&&\Sp_{12\ell}(2^f)\ar@{-}[dl]\ar@{}[d]|{\text{\cite[Theorem~A]{LPS1990}}}\ar@{-}[dr]&\\
&&&\Sp_{12}(2^{f\ell}){:}\ell\ar@{-}[dl]\ar@{-}[dr]&&\GO_{12\ell}^-(2^f)\ar@{-}[dl]\\
&&\Sp_{12}(2^{f\ell})\ar@{-}[dl]\ar@{-}[dr]&&\GO_{12}^-(2^{f\ell}){:}\ell\ar@{-}[dl]\\
&\Sp_6(2^{f\ell})\wr\Sy_2\ar@{-}[dl]\ar@{-}[dr]&&\GO_{12}^-(2^{f\ell}){:}\ell\ar@{-}[dl]&\\
G_2(2^{f\ell})\wr\Sy_2&&\GO_6^-(2^{f\ell})\times\GO_6^+(2^{f\ell})&&
}
}
$$
}
This gives a factorization $G=HK$ with $H=G_2(2^{f\ell})\wr\Sy_2=(\G_2(2^f)\times\G_2(2^f)){:}2$.

\begin{lemma}\label{lem17}
Let $W$ be a $2$-dimensional symplectic space over $\mathbb{F}_q$ with the non-degenerate symplectic form $\alpha$ and a standard basis $w_1,w_2$, and $V$ be a $2m$-dimensional symplectic space over $\mathbb{F}_q$ with the non-degenerate symplectic form $\beta$ and a standard basis $u_1,v_1,\dots,u_m,v_m$. Define a quadratic form $Q$ on $W\otimes V$ with the associated bilinear form:
$$
(W\otimes V)\times(W\otimes V)\rightarrow\mathbb{F}_q,\quad(w\otimes v,w'\otimes v')\mapsto\alpha(w,w')\beta(v,v')
$$
such that $Q(w\otimes v)=0$ for all $w\in W$ and $v\in V$. Let $B$ be the stabilizer in $\POm(W\otimes V,Q)$ of $\langle w_1\otimes u_1+w_2\otimes v_1\rangle$, and $M=\PSp_{2m}(q)$ be a subgroup of $\POm(W\otimes V,Q)$ that stabilizes $W\otimes V$ and acts trivially on $W$. Then, for $m\geqslant2$ and $q\geqslant4$,
$$
M\cap B=\{h\in\PSp_{2\ell}(q)\mid u_1^h=\lambda u_1,\ v_1^h=\lambda v_1,\ \lambda^2=1,\ \lambda\in\mathbb{F}_q^*\}.
$$
\end{lemma}

\begin{proof}
Direct calculation shows that
\begin{eqnarray*}
M\cap B&=&\{h\in\PSp_{2m}(q)\mid w_1\otimes u_1^h+w_2\otimes v_1^h=\lambda(w_1\otimes u_1+w_2\otimes v_1),\ \lambda\in\mathbb{F}_q^*\}\\
&=&\{h\in\PSp_{2m}(q)\mid w_1\otimes(u_1^h-\lambda u_1)=w_2\otimes(\lambda v_1-v_1^h),\ \lambda\in\mathbb{F}_q^*\}\\
&=&\{h\in\PSp_{2m}(q)\mid u_1^h=\lambda u_1,\ v_1^h=\lambda v_1,\ \lambda\in\mathbb{F}_q^*\}\\
&=&\{h\in\PSp_{2m}(q)\mid u_1^h=\lambda u_1,\ v_1^h=\lambda v_1,\ \beta(u_1^h,v_1^h)=\beta(u_1,v_1),\ \lambda\in\mathbb{F}_q^*\}\\
&=&\{h\in\PSp_{2m}(q)\mid u_1^h=\lambda u_1,\ v_1^h=\lambda v_1,\ \lambda^2=1,\ \lambda\in\mathbb{F}_q^*\},
\end{eqnarray*}
completing the proof.
\end{proof}

\begin{example}\label{exa5}
Let $G=\POm_{4\ell}^+(q)$ with $\ell\geqslant2$ and $q\geqslant4$, and $K=\N_1[G]=\Sp_{4\ell-2}(q)$. Then for each $U\leqslant\PSp_2(q)$, there exists a subgroup $H=U\times\PSp_{2\ell}(q)$ of $G$ such that $G=HK$, as in row~1 of Table~\ref{tab5}.
\end{example}

In fact, let $W$, $V$, $Q$, $B$ and $M$ be as defined in Lemma~\ref{lem17} with $m=\ell$. Let $G=\POm(W\otimes V,Q)$ and $K=B$ (as there is only one conjugacy class of subgroups $\N_1[G]$ in $G$). Take $H=U\times M$ to be a subgroup of $G$ stabilizing $W\otimes V$ which acts on $W$ and $V$ as $U$ and $M=\PSp_{2\ell}(q)$, respectively. To prove $G=HK$, we only need to prove $G=MK$ since $M\leqslant H$. By Lemma~\ref{lem17}, we have
\begin{eqnarray*}
|M\cap K||G|&=&|\{h\in\PSp_{2\ell}(q)\mid u_1^h=\lambda u_1,\ v_1^h=\lambda v_1,\ \lambda^2=1,\ \lambda\in\mathbb{F}_q^*\}||G|\\
&=&|\Sp_{2\ell-2}(q)||G|=|\PSp_{2\ell}(q)||K|=|M||K|.
\end{eqnarray*}
This shows that $G=MK$, proving Example~\ref{exa5}.

By~\cite[Theorem~A]{LPS1990}, there is a maximal factorization $G=AK$ with $G=\GaSp_{4\ell}(4)$, $A=\GaO_{4\ell}^+(4)$ and $K=\GaO_{4\ell}^-(4)$. We generalize this in the next lemma by showing that $A$ can be replaced by a smaller group $\Omega_{4\ell}^+(4).2$.

\begin{lemma}\label{lem18}
Let $G=L\rtimes\langle\phi\rangle=\GaSp_{4\ell}(4)$ with $\ell\geqslant2$, where $L=\Sp_{4\ell}(4)$ and $\phi$ is a field automorphism of order $2$, and $K=\GaO_{4\ell}^-(4)<G$. Then there exists $M=\Omega_{4\ell}^+(4)<L$ such that $G=(M\rtimes\langle\phi\rangle)K$ and $(M\rtimes\langle\phi\rangle)\cap K=\N_1[M]$.
\end{lemma}

\begin{proof}
From~\cite[3.2.4(e)]{LPS1990} we see that there is a subgroup $A=\GaO_{4\ell}^+(4)$ of $G$ such that $G=AK$, $A=(A\cap L)\rtimes\langle\phi\rangle$ and $A\cap K=\N_1[A\cap L]=\Sp_{4\ell-2}(4)\times2$. Taking $M=\Soc(A\cap L)$, we have $M=\Omega_{4\ell}^+(4)<L$, $A\cap L=\GO_{4\ell}^+(4)$ and $A=(M\rtimes\langle\phi\rangle)\N_1[A\cap L]$. This implies that $A=(M\rtimes\langle\phi\rangle)(A\cap K)$, and so $G=(M\rtimes\langle\phi\rangle)K$ by Lemma~\ref{lem16}. Moreover, $(M\rtimes\langle\phi\rangle)\cap K=(M\rtimes\langle\phi\rangle)\cap(A\cap K)$ since $M\rtimes\langle\phi\rangle<A$. Hence $(M\rtimes\langle\phi\rangle)\cap K=(M\rtimes\langle\phi\rangle)\cap\N_1[(A\cap L)]=\N_1[M]$, completing the proof.
\end{proof}

\begin{example}\label{exa2}
Let $G=\GaSp_{4\ell}(4)$ with $\ell\geqslant2$, and $K=\GaO_{4\ell}^-(4)<G$. Then there exists a subgroup $H=(\Sp_2(4)\times\Sp_{2\ell}(4)){:}2$ of $G$ such that $G=HK$, as in row~4 of Table~\ref{tab1}.
\end{example}

In fact, let $\phi$ and $M$ be as in Lemma~\ref{lem18} with $G=(M\rtimes\langle\phi\rangle)K$ and $(M\rtimes\langle\phi\rangle)\cap K=\N_1[M]$. As $M=\Omega_{4\ell}^+(4)$, there is a maximal subgroup $H=(\Sp_2(4)\otimes\Sp_{2\ell}(4))\rtimes\langle\phi\rangle$ of $M\rtimes\langle\phi\rangle$. Then $H=(\Sp_2(4)\times\Sp_{2\ell}(4)){:}2$, $M\rtimes\langle\phi\rangle=HM$, and $M=(H\cap M)\N_1[M]$ by~\cite[Theorem~A]{LPS1990}. It follows that $M\rtimes\langle\phi\rangle=H(H\cap M)\N_1[M]=H\N_1[M]$. Hence we have $G=HK$ as shown in the diagram
{\large
$$
\scalebox{0.8}[0.8]{
\xymatrix{
&&\GaSp_{4\ell}(4)\ar@{-}[dl]\ar@{-}[dr]&\\
&\Omega_{4\ell}^+(4)\rtimes\langle\phi\rangle\ar@{-}[dl]\ar@{-}[dr]&&\GaO_{4\ell}^-(4)\ar@{-}[dl]\\
(\Sp_2(4)\otimes\Sp_{2\ell}(4))\rtimes\langle\phi\rangle&&\N_1[\Omega_{4\ell}^+(4)]
}
}
$$
}

\begin{example}\label{exa3}
Let $G=\Sp_{8\ell}(2)$ and $K=\GO_{8\ell}^-(2)<G$. Then there exists a subgroup $H=(\Sp_2(4)\times\Sp_{2\ell}(4)){:}2$ of $G$ such that $G=HK$, as in row~5 of Table~\ref{tab1}. This is shown in the diagram
{\large
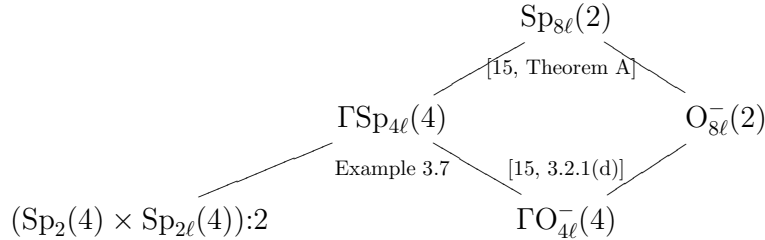
$$
\scalebox{0.8}[0.8]{
\xymatrix{
&&\Sp_{8\ell}(2)\ar@{-}[dl]\ar@{}[d]|{\text{\cite[Theorem~A]{LPS1990}}\ }\ar@{-}[dr]&\\
&\GaSp_{4\ell}(4)\ar@{}[d]|{\text{Example~\ref{exa2}}}\ar@{-}[dl]\ar@{-}[dr]
&\ar@{}[d]|{\text{\cite[3.2.1(d)]{LPS1990}}}&\GO_{8\ell}^-(2)\ar@{-}[dl]\\
(\Sp_2(4)\times\Sp_{2\ell}(4)){:}2&&\GaO_{4\ell}^-(4)
}
}
$$
}
\end{example}

Before stating the next lemma, recall from~\cite[Theorem~A]{LPS1990} that there is a maximal factorization $G=AK$ with $G=\Omega_{4m}^+(2)$, $A=\GaO_{2m}^+(4)$ and $K=\N_1[G]$. We generalize this by showing in the lemma below that $A$ can be replaced by a smaller group $\Omega_{2m}^+(4).2$.

\begin{lemma}\label{lem19}
Let $G=\Omega_{4m}^+(2)$ with $m\geqslant2$, $K=\N_1[G]=\Sp_{4m-2}(2)$, $A=\GaO_{2m}^+(4)<G$ and $M=\Soc(A)\rtimes\langle\phi\rangle<G$, where $\phi$ is the field automorphism of $\mathbb{F}_4$ of order $2$. Then $A\cap K\leqslant\GO_{2m}^+(4)$, $G=MK$ and $M\cap K=\N_1[\Soc(A)]=\Sp_{2m-2}(4)$.
\end{lemma}

\begin{proof}
Let $V$ be a $2m$-dimensional orthogonal space over $\mathbb{F}_4$ with a non-degenerate quadratic form $P$ of plus type, and let $e_1,f_1,\dots,e_m,f_m$ be a standard basis of $V$. Take $A=\GO(V,P)\rtimes\langle\phi\rangle$, so that $M=\Omega(V,P)\rtimes\langle\phi\rangle$ and $\phi$ fixes each of $e_1,f_1,\dots,e_m,f_m$. Define a map $Q$ from $V$ to $\mathbb{F}_2$ by putting $Q(u)=P(u)+P(u)^2$ for each $u\in V$. Then $Q$ is a non-degenerate quadratic form of plus type on $V$ while $V$ is considered as a $4m$-dimensional orthogonal space over $\mathbb{F}_2$ with respect to $Q$. Regard $G$ as $\Omega(V,Q)$ and regard $K$ as the stabilizer in $G$ of $e_1+\omega f_1$, where $\omega$ is a generator of $\mathbb{F}_4^\times$. (Note that $Q(e_1+\omega f_1)=P(e_1+\omega f_1)+P(e_1+\omega f_1)^2=\omega+\omega^2=1$.) Then $A<G$ since both $\GO(V,P)$ and $\phi$ preserve $P(u)+P(u)^2$ for each $u\in V$.

We first prove $A\cap K\leqslant\GO(V,P)$. To see this, suppose on the contrary that there exists $g\in\GO(V,P)\setminus(A\cap K)$. Then since $A=\GO(V,P)\rtimes\langle\phi\rangle$, we conclude that $g\phi\in K$. Writing $(e_1+\omega f_1)^g=\sum_{i=1}^{2\ell}(a_ie_i+b_if_i)$ with $a_1,b_1,\dots,a_{2\ell},b_{2\ell}\in\mathbb{F}_4$, we deduce that
$$
e_1+\omega f_1=(e_1+\omega f_1)^{g\phi}=\left(\sum_{i=1}^{2\ell}(a_ie_i+b_if_i)\right)^\phi=\sum_{i=1}^{2\ell}(a_i^2e_i+b_i^2f_i).
$$
This means $a_1^2=1$, $b_1^2=\omega$ and $a_i^2=b_i^2=0$ for $2\leqslant i\leqslant2\ell$, which is equivalent to $a_1=1$, $b_1=\omega^2$ and $a_i=b_i=0$ for $2\leqslant i\leqslant2\ell$. Now
$$
(e_1+\omega f_1)^g=\sum_{i=1}^{2\ell}(a_ie_i+b_if_i)=e_1+\omega^2f_1,
$$
but $P(e_1+\omega f_1)=\omega\neq\omega^2=P(e_1+\omega^2f_1)$, contradicting $g\in\GO(V,P)$. Consequently, $A\cap K\leqslant\GO(V,P)=\GO_{2m}^+(4)$.

By the above conclusion we know that $M\cap K\leqslant\GO(V,P)$, and so
$$
M\cap K=\GO(V,P)\cap M\cap K=\Omega(V,P)\cap K=\N_1[\Soc(A)]=\Sp_{2m-2}(4).
$$
As a consequence, $|M\cap K||G|=|M||K|$. This yields $G=MK$ by Lemma~\ref{lem6}.
\end{proof}

\begin{example}\label{exa10}
Let $G=\Omega_{8\ell}^+(2)$ and $K=\N_1[G]=\Sp_{8\ell-2}(2)$. Then there exists a subgroup $H=(\Sp_2(4)\times\Sp_{2\ell}(4)){:}2$ of $G$ such that $G=HK$, as in row~2 of Table~\ref{tab5}.
\end{example}

In fact, let $\phi$ and $M$ be as in Lemma~\ref{lem19} with $m=2\ell$, so that $G=MK$ and $M\cap K=\N_1[\Soc(M)]=\Sp_{4\ell-2}(4)$. Taking $H=(\Sp_2(4)\otimes\Sp_{2\ell}(4))\rtimes\langle\phi\rangle$ to be a maximal subgroup of $M=\Omega_{4\ell}^+(4)\rtimes\langle\phi\rangle$, we have $H=(\Sp_2(4)\times\Sp_{2\ell}(4)){:}2$, $M=H\Soc(M)$, and $\Soc(M)=(H\cap\Soc(M))\N_1[\Soc(M)]$ by~\cite[Theorem~A]{LPS1990}. It follows that $M=H(H\cap\Soc(M))\N_1[\Soc(M)]=H\N_1[\Soc(M)]$, and then the diagram
{\large
$$
\scalebox{0.8}[0.8]{
\xymatrix{
&&\Omega_{8\ell}^+(2)\ar@{-}[dl]\ar@{-}[dr]&\\
&\Omega_{4\ell}^+(4)\rtimes\langle\phi\rangle\ar@{-}[dl]\ar@{-}[dr]&&\Sp_{8\ell-2}(2)\ar@{-}[dl]\\
(\Sp_2(4)\otimes\Sp_{2\ell}(4))\rtimes\langle\phi\rangle&&\N_1[\Omega_{4\ell}^+(4)]
}
}
$$
}
shows that $G=HK$.

\begin{example}\label{exa6}
Let $G=\GaO_{8\ell}^+(4)$ and $K=\N_1[G]$. Then there exists a subgroup $H=(\Sp_2(4)\times\Sp_{2\ell}(16)){:}4$ of $G$ such that $G=HK$, as in row~3 of Table~\ref{tab5}.
\end{example}

Let $W$, $V$, $Q$ and $M$ be as in Lemma~\ref{lem17} with $m=2\ell$ and $q=4$, and take
$$
e_i=w_1\otimes u_i,\quad f_i=w_2\otimes v_i,\quad e_{i+2\ell}=w_1\otimes v_i\quad\text{and}\quad f_{i+2\ell}=w_2\otimes u_i
$$
for $1\leqslant i\leqslant2\ell$. Then $e_1,f_1,\dots,e_{4\ell},f_{4\ell}$ is a standard basis of the orthogonal space $W\otimes V$. Let $\phi$ be the field automorphism of order $2$ fixing each of $e_1,f_1,\dots,e_{4\ell},f_{4\ell}$. Then we can regard $G$ as $\GO(W\otimes V,Q)\rtimes\langle\phi\rangle$ and regard $K$ as the stabilizer in $G$ of $\langle e_1+f_1\rangle$. Let $A=(C\times M)\rtimes\langle\phi\rangle$ be the subgroup of $G$ stabilizing $W\otimes V$, where $C=\Sp_2(4)$ acts on $V$ trivially, and denote the actions of $A$ on $W$ and $V$ by $\pi_1$ and $\pi_2$ respectively. It is evident that $A^{\pi_1}=\GaSp_2(4)$ and $A^{\pi_2}=\GaSp_{4\ell}(4)$. Take $D=\Sp_{2\ell}(16){:}2$ to be a $\calC_3$-subgroup of $M$ such that $(D\rtimes\langle\phi\rangle)^{\pi_2}=\GaSp_{2\ell}(16)$. Then there exists an element $\sigma$ of $M$ with $(D\rtimes\langle\phi\rangle)^{\pi_2}=\Soc(D)\rtimes\langle\sigma\phi\rangle^{\pi_2}=\Sp_{2\ell}(16){:}4$. As $(\sigma\phi)^{\pi_1}=\phi^{\pi_1}$ has order $2$ and $(\sigma\phi)^{\pi_2}$ has order $4$, we see that $\sigma\phi$ has order $4$. Also, $\langle\sigma\phi\rangle\cap\Soc(D)=1$ since $((\sigma\phi)^2)^{\pi_2}\notin\Soc(D)$. Furthermore, $\sigma\phi$ normalizes $C\times\Soc(D)$ since both $\sigma$ and $\phi$ normalize $C\times\Soc(D)$. Thus we have a subgroup $H=(C\times\Soc(D))\rtimes\langle\sigma\phi\rangle=(\Sp_2(4)\times\Sp_{2\ell}(16)){:}4$ of $A$.

We claim that $(A\cap K)^{\pi_2}$ contains the stabilizer in $A^{\pi_2}$ of the non-degenerate $2$-space $\langle u_1,v_1\rangle$. Since $A\cap K\geqslant\langle\phi\rangle$, it suffices to show that $(A\cap K)^{\pi_2}$ contains the stabilizer of $\langle u_1,v_1\rangle$ in $M^{\pi_2}=\Sp_{4\ell}(4)$. Let $g$ be an arbitrary element of $M^{\pi_2}$ that stabilizes $\langle u_1,v_1\rangle$. Then $u_1^g=au_1+bv_1$ and $v_1^g=cu_1+dv_1$ for some $a,b,c,d\in\mathbb{F}_4$ with $ad+bc=1$. Take $x\in C$ and $y\in M$ such that $w_1^x=dw_1+cw_2$, $w_2^x=bw_1+aw_2$ and $y^{\pi_2}=g$, and let $k=xy\in A$. It is straightforward to verify that $k^{\pi_2}=g$ and $k$ fixes $e_1+f_1$, which means $k\in A\cap K$. This proves our claim.

Now we have $H^{\pi_2}\geqslant\Soc(D)\rtimes\langle\sigma\phi\rangle^{\pi_2}=(D\rtimes\langle\phi\rangle)^{\pi_2}=\GaSp_{2\ell}(16)$ and $(A\cap K)^{\pi_2}\geqslant\N_2[A^{\pi_2}]$. It follows from~\cite[Theorem~A]{LPS1990} that $A^{\pi_2}=\GaSp_{4\ell}(4)$ has a factorization $A^{\pi_2}=H^{\pi_2}(A\cap K)^{\pi_2}$. Then since $A$ contains $C=\ker(\pi_2)$, we deduce that $A=H(A\cap K)$. Therefore, the diagram
{\large
$$
\scalebox{0.8}[0.8]{
\xymatrix{
&&G=\GaO_{8\ell}^+(4)\ar@{-}[dl]\ar@{}[d]|{\text{\cite[Theorem~A]{LPS1990}}\ }\ar@{-}[dr]&\\
&A=(\Sp_2(4)\times\Sp_{4\ell}(4)){:}2\ar@{-}[dl]\ar@{-}[dr]&&K=\N_1\ar@{-}[dl]\\
H&&A\cap K
}
}
$$
}
shows that $G=HK$, proving Example~\ref{exa6}.

\begin{example}\label{exa4}
Let $L=\Sp_6(2^f)$ with $f\geqslant2$, $G$ be an almost simple group with socle $L$, $H=U\times\Sp_4(2^f)\leqslant\N_2[L]$ with $U\leqslant\Sp_2(2^f)$, and $K$ be a maximal subgroup of $G$ such that $K\cap L=\G_2(2^f)$. Then $G=HK$, as in row~7 of Table~\ref{tab1}.
\end{example}

By~\cite[4.3]{Wilson2009}, there is a $6$-dimensional vector space $V$ over $\mathbb{F}_{2^f}$ with basis
$$
\overline{x_1},\overline{x_2},\overline{x_3},\overline{x_6},\overline{x_7},\overline{x_8}
$$
and a bilinear form $\beta$ on $V$ with Gram matrix
$$
\begin{pmatrix}
&&&&&1\\
&&&&1&\\
&&&1&&\\
&&1&&&\\
&1&&&&\\
1&&&&&
\end{pmatrix}
$$
such that $L:=\Sp(V,\beta)=\Sp_6(2^f)$ has a subgroup $N=\G_2(2^f)$ generated by
\begin{eqnarray*}
&r:\ (\overline{x_1},\overline{x_2},\overline{x_3},\overline{x_6},\overline{x_7},\overline{x_8})\mapsto
(\overline{x_1},\overline{x_3},\overline{x_2},\overline{x_7},\overline{x_6},\overline{x_8}),\\
&s:\ (\overline{x_1},\overline{x_2},\overline{x_3},\overline{x_6},\overline{x_7},\overline{x_8})\mapsto
(\overline{x_2},\overline{x_1},\overline{x_6},\overline{x_3},\overline{x_8},\overline{x_7}),\\
&T(\lambda,\mu):\ (\overline{x_1},\overline{x_2},\overline{x_3},\overline{x_6},\overline{x_7},\overline{x_8})\mapsto
(\lambda\overline{x_1},\mu\overline{x_2},\lambda\mu^{-1}\overline{x_3},
\lambda^{-1}\mu\overline{x_6},\mu^{-1}\overline{x_7},\lambda^{-1}\overline{x_8}),\\
&A(\alpha):\ \overline{x_7}\mapsto\overline{x_7}+\alpha\overline{x_1},\quad\overline{x_8}\mapsto\overline{x_8}+\alpha\overline{x_2}
\end{eqnarray*}
and
$$
F(\alpha):\ \overline{x_2}\mapsto\overline{x_2}+\alpha\overline{x_1},
\quad\overline{x_6}\mapsto\overline{x_6}+\alpha^2\overline{x_3},\quad\overline{x_8}\mapsto\overline{x_8}+\alpha\overline{x_7},
$$
where $\lambda$ runs over $\mathbb{F}_{2^f}^*$, $\mu$ runs over $\mathbb{F}_{2^f}^*$ and $\alpha$ runs over $\mathbb{F}_{2^f}$. Let $M$ be the subgroup of $L$ generated by $s,(rs)^3,T(\lambda,\mu),A(\alpha),F(\alpha)$, and $\phi$ be the semilinear transformation of $V$ defined by
$$
(\alpha_1\overline{x_1}+\alpha_2\overline{x_2}+\alpha_3\overline{x_3}+\alpha_6\overline{x_6}+\alpha_7\overline{x_7}+\alpha_8\overline{x_8})^\phi
=\alpha_1^2\overline{x_1}+\alpha_2^2\overline{x_2}+\alpha_3^2\overline{x_3}
+\alpha_6^2\overline{x_6}+\alpha_7^2\overline{x_7}+\alpha_8^2\overline{x_8}
$$
for $\alpha_1,\alpha_2,\alpha_3,\alpha_6,\alpha_7,\alpha_8\in\mathbb{F}_{2^f}$. Then we can write $G=L\rtimes\langle\phi^e\rangle$ and $K=N\rtimes\langle\phi^e\rangle$ for some divisor $e$ of $f$. Moreover, $M$ stabilizes $\langle\overline{x_3},\overline{x_6}\rangle$ and $\langle\overline{x_1},\overline{x_2},\overline{x_7},\overline{x_8}\rangle$, respectively. Let $A$ be the stabilizer in $G$ of $\langle\overline{x_3},\overline{x_6}\rangle$ whose actions on $\langle\overline{x_3},\overline{x_6}\rangle$ and $\langle\overline{x_1},\overline{x_2},\overline{x_7},\overline{x_8}\rangle$ are denoted by $\pi_1$ and $\pi_2$, respectively. It follows that $K\cap L=\G_2(2^f)$, and $A\cap L=(A\cap L)^{\pi_1}\times(A\cap L)^{\pi_2}$ with $(A\cap L)^{\pi_1}=\Sp_2(2^f)$ and $(A\cap L)^{\pi_2}=\Sp_4(2^f)$. Note that $M\leqslant A\cap K\cap L$ while $M=\Omega_4^+(2^f)$ is a maximal subgroup of $K\cap L$. We conclude that $A\cap K\cap L=M$ since $K\cap L$ does not stabilize $\langle\overline{x_3},\overline{x_6}\rangle$. Now as
$$
|A\cap K\cap L||L|=|M||L|=|A\cap L||K\cap L|,
$$
we have $L=(A\cap L)(K\cap L)$ by Lemma~\ref{lem6}. Regard $H=U\times(A\cap L)^{\pi_2}\leqslant A\cap L$, $U$ being a subgroup of $(A\cap L)^{\pi_1}$.

Note that the subgroup of $M$ generated by $(rs)^3,T(\lambda,1),F(\alpha)$, where $\lambda$ runs over $\mathbb{F}_{2^f}^*$ and $\alpha$ runs over $\mathbb{F}_{2^f}$, acts on $\langle\overline{x_3},\overline{x_6}\rangle$ as $\Sp_2(2^f)$. We have $M^{\pi_2}\geqslant(A\cap L)^{\pi_1}$ and so $(A\cap K\cap L)^{\pi_1}=(A\cap L)^{\pi_1}$. Accordingly,
\begin{eqnarray*}
H(A\cap K\cap L)&\geqslant&(A\cap L)^{\pi_2}(A\cap K\cap L)\\
&=&(A\cap L)^{\pi_2}(A\cap K\cap L)^{\pi_1}=(A\cap L)^{\pi_2}(A\cap L)^{\pi_1}=A\cap L,
\end{eqnarray*}
which means $A\cap L=H(A\cap K\cap L)$. This together with $H\leqslant A\cap L$ and $L=(A\cap L)(K\cap L)$ implies that $L=H(K\cap L)$. Hence $HK=H(K\cap L)\langle\phi^e\rangle=L\langle\phi^e\rangle=G$, which proves Example~\ref{exa4}.

\begin{example}\label{exa11}
Let $G=\Omega_{12}^+(2^f)$ with $f\geqslant2$, and $K=\N_1[G]=\Sp_{10}(2^f)$. Then for each $U\leqslant\Sp_2(2^f)$, there exists a subgroup $H=U\times\G_2(2^f)$ of $G$ such that $G=HK$, as in row~4 of Table~\ref{tab5}.
\end{example}

Let $W$, $V$, $Q$, $B$ and $M$ be as defined in Lemma~\ref{lem17} with $m=3$ and $q=2^f$. Then we may let $G=\Omega(W\otimes V,Q)$ and $K=B$. It follows from Lemma~\ref{lem17} that $M\cap K=\Sp_4(2^f)<\N_2[M]$. Hence $|M\cap K||G|=|M||K|$, and so $G=MK$ by Lemma~\ref{lem6}. Let $N=\G_2(2^f)<M$, and $H=U\otimes N$ be a subgroup of $G$ stabilizing $W\otimes V$ which acts on $W$ and $V$ as $U$ and $N$, respectively. To prove that $G=HK$, it suffices to prove $G=NK$ since $N\leqslant H$. In fact, Example~\ref{exa4} shows that $M=N(M\cap K)$, which together with $G=MK$ leads to $G=NK$ by Lemma~\ref{lem16}. This implies $G=HK$, proving Example~\ref{exa11}. We remark that here $M$ is contained in a maximal subgroup $A=\Sp_2(2^f)\otimes\Sp_6(2^f)$ of $G$, so that we have the diagram
{\large
$$
\scalebox{0.8}[0.8]{
\xymatrix{
&&&G\ar@{-}[dl]\ar@{}[d]\ar@{-}[dr]&\\
&&A\ar@{-}[dl]\ar@{-}[dr]&&K\ar@{-}[dl]\\
&M\ar@{-}[dl]\ar@{-}[dr]&&A\cap K\ar@{-}[dl]\\
H&&M\cap K&
}
}
$$
}

\begin{example}\label{exa12}
Let $G=\GaSp_{12}(4)$ and $H=\N_2[G]=(\Sp_2(4)\times\Sp_{10}(4)){:}2$. Then there exists a subgroup $K=\G_2(16){:}4$ of $G$ such that $G=HK$, as in row~8 of Table~\ref{tab1}.
\end{example}

Let $B=\GaSp_6(16)$ be a $\calC_3$-subgroup of $G$. According to~\cite[3.2.1(a)]{LPS1990} we have $G=HB$ with $H\cap B=H\cap\Soc(B)=\Sp_2(4)\times\Sp_4(16)$. Take $K=\G_2(16){:}4$ to be a maximal subgroup of $B$ (see~\cite[Table~8.29]{BHR2013}). Then $B=\Soc(B)K$ and $K\cap\Soc(B)=\G_2(16)$. Moreover, $\Soc(B)=(H\cap\Soc(B))(K\cap B)$ as Example~\ref{exa4} shows. Thereby we deduce $B=(H\cap\Soc(B))K$ and thus $B=(H\cap B)K$. This implies $G=HK$ as shown in the diagram
{\large
$$
\scalebox{0.8}[0.8]{
\xymatrix{
&G=\GaSp_{12}(4)\ar@{-}[dl]\ar@{}[d]\ar@{-}[dr]&&&\\
H=\N_2\ar@{-}[dr]&&B=\GaSp_6(16)\ar@{-}[dl]\ar@{-}[dr]&\\
&H\cap B=\Sp_2(4)\times\Sp_4(16)&&K=\G_2(16){:}4&
}
}
$$
}
proving Example~\ref{exa12}.

\begin{example}\label{exa13}
Let $G=\GaSp_{12}(4)$ and $K=\GaO_{12}^-(4)<G$. Then there exists a subgroup $H=(\Sp_2(4)\times\G_2(4)){:}2$ of $G$ such that $G=HK$, as in row~9 of Table~\ref{tab1}.
\end{example}

Let $\phi$ and $M$ be as in Lemma~\ref{lem18} with $\ell=3$, so that $G=(M\rtimes\langle\phi\rangle)K$ and $(M\rtimes\langle\phi\rangle)\cap K=\N_1[M]$. Following the proof of Example~\ref{exa11} we obtain a subgroup $H^*=\Sp_2(4)\times\G_2(4)$ of $M$ such that $M=H^*\N_1[M]$ and $H^*$ is normalized by $\phi$. Take $H=H^*\rtimes\langle\phi\rangle<M\rtimes\langle\phi\rangle$. We have $H=(\Sp_2(4)\times\G_2(4)){:}2$ and $M\rtimes\langle\phi\rangle=HM$. It follows from Lemma~\ref{lem16} that $M\rtimes\langle\phi\rangle=H\N_1[M]$ since $H\cap M=H^*$ and $M=H^*\N_1[M]$. This implies $G=HK$ as shown in the diagram
{\large
$$
\scalebox{0.8}[0.8]{
\xymatrix{
&&G=\GaSp_{12}(4)\ar@{-}[dl]\ar@{-}[dr]&\\
&M\rtimes\langle\phi\rangle=\Omega_{12}^+(4){:}2\ar@{-}[dl]\ar@{-}[dr]&&K=\GaO_{12}^-(4)\ar@{-}[dl]\\
H&&\N_1[M]
}
}
$$
}
proving Example~\ref{exa13}.

\begin{example}\label{exa14}
Let $G=\Omega_{24}^+(2)$ and $K=\N_1[G]=\Sp_{22}(2)$. Then there exists a subgroup $H=(\Sp_2(4)\times\G_2(4)){:}2$ of $G$ such that $G=HK$, as in row~5 of Table~\ref{tab5}.
\end{example}

Let $\phi$ and $M$ be as in Lemma~\ref{lem19} with $m=6$, so that $G=MK$ and $M\cap K=\N_1[\Soc(M)]$. Following the proof of Example~\ref{exa11} we obtain a subgroup $H^*=\Sp_2(4)\times\G_2(4)$ of $\Soc(M)$ such that $\Soc(M)=H^*\N_1[\Soc(M)]$ and $H^*$ is normalized by $\phi$. Take $H=H^*\rtimes\langle\phi\rangle<M$. We have $H=(\Sp_2(4)\times\G_2(4)){:}2$ and $M=H\Soc(M)$. This implies $M=H\N_1[\Soc(M)]$ by Lemma~\ref{lem16} since $H\cap\Soc(M)=H^*$ and $\Soc(M)=H^*\N_1[\Soc(M)]$. Then the diagram 
{\large
$$
\scalebox{0.8}[0.8]{
\xymatrix{
&&G=\Omega_{24}^+(2)\ar@{-}[dl]\ar@{-}[dr]&\\
&M=\Omega_{12}^+(4){:}2\ar@{-}[dl]\ar@{-}[dr]&&K=\N_1[G]\ar@{-}[dl]\\
H&&\N_1[\Soc(M)]
}
}
$$
}
shows that $G=HK$, proving Example~\ref{exa14}.

\begin{example}\label{exa15}
Let $G=\Sp_{24}(2)$ and $K=\GO_{24}^-(2)<G$. Then there exists a subgroup $H=(\Sp_2(4)\times\G_2(4)){:}2$ of $G$ such that $G=HK$, as in row~10 of Table~\ref{tab1}. This is shown in the diagram
{\large
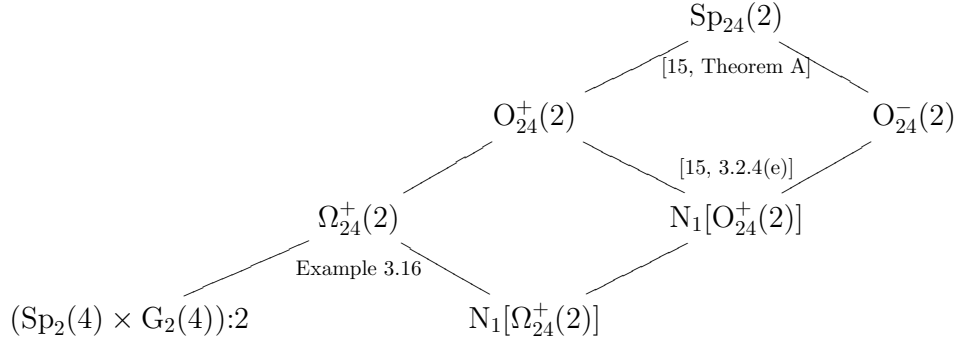
$$
\scalebox{0.8}[0.8]{
\xymatrix{
&&&\Sp_{24}(2)\ar@{-}[dl]\ar@{}[d]|{\text{\cite[Theorem~A]{LPS1990}}}\ar@{-}[dr]&\\
&&\GO_{24}^+(2)\ar@{-}[dl]\ar@{-}[dr]&\ar@{}[d]|{\text{\cite[3.2.4(e)]{LPS1990}}}&\GO_{24}^-(2)\ar@{-}[dl]\\
&\Omega_{24}^+(2)\ar@{}[d]|{\text{Example~\ref{exa14}}}\ar@{-}[dl]\ar@{-}[dr]&&\N_1[\GO_{24}^+(2)]\ar@{-}[dl]\\
(\Sp_2(4)\times\G_2(4)){:}2&&\N_1[\Omega_{24}^+(2)]&
}
}
$$
}
\end{example}

\begin{example}\label{exa16}
Let $G=\GaO_{24}^+(4)$ and $K=\N_1[G]$. Then there exists a subgroup $H=(\Sp_2(4)\times\G_2(16)){:}4$ of $G$ such that $G=HK$, as in row~6 of Table~\ref{tab5}.
\end{example}

Let $W$, $V$, $Q$ and $M$ be as in Lemma~\ref{lem17} with $m=6$ and $q=4$, and take
$$
e_i=w_1\otimes u_i,\quad f_i=w_2\otimes v_i,\quad e_{i+6}=w_1\otimes v_i\quad\text{and}\quad f_{i+6}=w_2\otimes u_i
$$
for $1\leqslant i\leqslant6$. Then $e_1,f_1,\dots,e_{12},f_{12}$ is a standard basis of the orthogonal space $W\otimes V$. Let $\phi$ be the field automorphism of order $2$ fixing each of $e_1,f_1,\dots,e_{12},f_{12}$. Then we can regard $G$ as $\GO(W\otimes V,Q)\rtimes\langle\phi\rangle$ and regard $K$ as the stabilizer in $G$ of $\langle e_1+f_1\rangle$. Let $A=(C\times M)\rtimes\langle\phi\rangle$ be the subgroup of $G$ stabilizing $W\otimes V$, where $C=\Sp_2(4)$ acts on $V$ trivially, and denote the actions of $A$ on $W$ and $V$ by $\pi_1$ and $\pi_2$ respectively. It is evident that $A^{\pi_1}=\GaSp_2(4)$ and $A^{\pi_2}=\GaSp_{12}(4)$. Following the proof of Example~\ref{exa12} we see that there is a subgroup $D=\G_2(16){:}2$ of $M=\Sp_{12}(4)$ with $(D\rtimes\langle\phi\rangle)^{\pi_2}=\G_2(16){:}4$ and $A^{\pi_2}=(D\rtimes\langle\phi\rangle)^{\pi_2}\N_2[A^{\pi_2}]$. Accordingly, there exists an element $\sigma$ of $M$ such that $(D\rtimes\langle\phi\rangle)^{\pi_2}=\Soc(D)\rtimes\langle\sigma\phi\rangle^{\pi_2}=\G_2(16){:}4$. Since $(\sigma\phi)^{\pi_1}=\phi^{\pi_1}$ has order $2$ and $(\sigma\phi)^{\pi_2}$ has order $4$, we conclude that $\sigma\phi$ has order $4$. Also, $\langle\sigma\phi\rangle\cap\Soc(D)=1$ as $((\sigma\phi)^2)^{\pi_2}\notin\Soc(D)$. Furthermore, $\sigma\phi$ normalizes $C\times\Soc(D)$ since both $\sigma$ and $\phi$ normalize $C\times\Soc(D)$. Thus we have a subgroup $H=(C\times\Soc(D))\rtimes\langle\sigma\phi\rangle=(\Sp_2(4)\times\G_2(16)){:}4$ of $A$.

As in the proof of Example~\ref{exa6} we know that $(A\cap K)^{\pi_2}$ contains the stabilizer in $A^{\pi_2}$ of the non-degenerate $2$-space $\langle u_1,v_1\rangle$. Now $(A\cap K)^{\pi_2}\geqslant\N_2(A^{\pi_2})$ and $A^{\pi_2}=H^{\pi_2}\N_2[A^{\pi_2}]$. It follows that $A^{\pi_2}=H^{\pi_2}(A\cap K)^{\pi_2}$, and so $A=H(A\cap K)$ since $A$ contains $C=\ker(\pi_2)$. Therefore, the diagram
{\large
$$
\scalebox{0.8}[0.8]{
\xymatrix{
&&G=\GaO_{24}^+(4)\ar@{-}[dl]\ar@{}[d]|{\text{\cite[Theorem~A]{LPS1990}}\ }\ar@{-}[dr]&\\
&A=(\Sp_2(4)\times\Sp_{12}(4)){:}2\ar@{-}[dl]\ar@{-}[dr]&&K=\N_1\ar@{-}[dl]\\
H&&A\cap K
}
}
$$
}
shows that $G=HK$, proving Example~\ref{exa16}.

\begin{corollary}
For each $L$ as in the second column of \emph{Table~\ref{tab1}}, there exist an almost simple group $G$ with socle $L$ and a factorization $G=HK$ such that $H\cap L$ and $K\cap L$ are as described in the same row as $L$ in \emph{Table~\ref{tab1}}.
\end{corollary}

\begin{proof}
For row~3 and row~6 of Table~\ref{tab1}, there is a maximal factorization $G=HK$ (corresponding to row~2 and row~3 of Table~\ref{tab2} respectively) as Lemma~\ref{lem1} shows. For the other rows of Table~\ref{tab1}, the assertion is true by the example indicated in the last column of Table~\ref{tab1}.
\end{proof}

\section{Strategy of proof}

As mentioned in the Introduction, the maximal factorizations of almost simple groups have been classified by Liebeck, Praeger and Saxl~\cite{LPS1990}. In order to apply this result to investigate the general factorizations of an almost simple group $G$, say, we need to embed a given factorization $G=HK$ to the a maximal factorization $G=AB$. This may be easily accomplished by taking arbitrary maximal subgroups $A$ and $B$ of $G$ containing $H$ and $K$ respectively. However, such maximal subgroups $A$ and $B$ are not necessarily core-free even if $H$ and $K$ are core-free. For example, if $HL<G$, where $L=\Soc(G)$, then the maximal subgroup of $G$ containing $HL$ (and thus containing $H$) is not core-free in $G$. In fact, $HL=G$ if and only if all maximal subgroups of $G$ containing $H$ are core-free in $G$, as proved in the following lemma.

\begin{lemma}\label{Maximality}
Let $G$ be an almost simple group with socle $L$ and $H$ be a core-free subgroups of $G$. Then $HL=G$ if and only if all maximal subgroups of $G$ containing $H$ are core-free in $G$.
\end{lemma}

\begin{proof}
Suppose that $HL=G$. For any maximal subgroup $A$ of $G$ containing $H$, since $A\geqslant L$ would lead to a contradiction that $A\geqslant HL=G$, we see that $A$ is core-free in $G$.

Conversely, suppose that all maximal subgroups of $G$ containing $H$ are core-free in $G$. If $HL<G$, then the maximal subgroup of $G$ containing $HL$ (and thus containing $H$) is not core-free in $G$, a contradiction. Hence $HL=G$. This proves the lemma.
\end{proof}

The next lemma is a direct application of \cite[Lemma~2(i)]{LPS1996}.

\begin{lemma}\label{Embedding}
Suppose that $G$ is an almost simple group with socle $L$ and $G=HK$ with core-free subgroups $H$ and $K$. Then there exists a factorization $G^*=H^*K^*$ such that $L\trianglelefteq G^*\leqslant G$, $H^*\cap L=H\cap L$, $K^*\cap L=K\cap L$ and $H^*L=K^*L=G^*$.
\end{lemma}

\begin{proof}
Take $G^*=HL\cap KL$, $H^*=H\cap G^*$ and $K^*=K\cap G^*$. Then $L\trianglelefteq G^*\leqslant G$, $H^*\cap L=H\cap L$ and $K^*\cap L=K\cap L$. It follows that $H^*$ and $K^*$ are both core-free in $G^*$ since $H$ and $K$ are both core-free in $G$. By \cite[Lemma~2(i)]{LPS1996}, we have $G^*=H^*K^*$ and $H^*L=K^*L=G^*$. Thus, $G^*=H^*K^*$ is a factorization as desired.
\end{proof}

To describe $H\cap\Soc(G)$ and $K\cap\Soc(G)$ (as Theorem~\ref{thm1}(b) does) for a factorization $G=HK$ of an almost simple group $G$ with $H$ and $K$ core-free, we may assume by virtue of Lemma~\ref{Embedding} that $HL=KL=G$. According to Lemma~\ref{Maximality}, this is equivalent to assuming that all maximal subgroups of $G$ containing $H$ and $K$, respectively, are core-free in $G$. After we know $H\cap\Soc(G)$ and $K\cap\Soc(G)$, the following lemma will help to determine the factors $H$ and $K$ of the factorization $G=HK$.

\begin{lemma}
Suppose that $G$ is an almost simple group with socle $L$ and $G=HK$ with core-free subgroups $H$ and $K$. Then there exists a factorization $G/L=\overline{H}\,\overline{K}$ such that $H=(H\cap L).\overline{H}$ and $K=(K\cap L).\overline{K}$.
\end{lemma}

\begin{proof}
Let $\overline{H}=HL/L$ and $\overline{K}=KL/L$. Then by modulo $L$, we derive from $G=HK$ that $G/L=\overline{H}\,\overline{K}$. Moreover, since $\overline{H}=HL/L\cong H/(H\cap L)$ and $\overline{K}=KL/L\cong K/(K\cap L)$, we have $H=(H\cap L).\overline{H}$ and $K=(K\cap L).\overline{K}$. Hence the lemma follows.
\end{proof}

\begin{remark}
In an alternative way, one may embed a factorization $G=HK$ of an almost simple group $G$ with $H$ and $K$ core-free into a factorization of $G$ with both factors maximal among core-free subgroups of $G$. The latter is called a $\mathrm{max}^-$ factorization of $G$ and is characterized in~\cite{LPS1996}.
\end{remark}

Now let $G$ be an almost simple group with socle $L$ and let $G=HK$ be a factorization of $G$ such that $H$ has at least two nonsolvable composition factors and $K$ is core-free, as in the assumption of Theorem~\ref{thm1}. We first prove the conclusion of Theorem~\ref{thm1} in the case that $L$ is a non-classical simple group.

\begin{lemma}\label{lem2}
If $L$ is a non-classical simple group, then $L=\A_n$ with $n\geqslant10$ and one of the following holds:
\begin{itemize}
\item[(a)] $H$ is a transitive permutation group of degree $n$, and $\A_{n-1}\leqslant K\leqslant\Sy_{n-1}$;
\item[(b)] $n=10$, $H$ is a transitive permutation group of degree $10$ such that $(\A_5\times\A_5).2\leqslant H\leqslant\Sy_5\wr\Sy_2$, and $K=\SL_2(8)$ or $\SL_2(8).3$;
\item[(c)] $n=12$, $\A_7\times\A_5\leqslant H\leqslant\Sy_7\times\Sy_5$, and $K=\M_{12}$;
\item[(d)] $n=24$, $\A_{19}\times\A_5\leqslant H\leqslant\Sy_{19}\times\Sy_5$, and $K=\M_{24}$.
\end{itemize}
\end{lemma}

\begin{proof}
Let $L$ be a non-classical simple group satisfying the assumption of Theorem~\ref{thm1}. Consulting the classification of factorizations of exceptional groups of Lie type in~\cite{HLS1987}, one sees that $L$ is not an exceptional group of Lie type as $H$ has at least two nonsolvable composition factors. Similarly, $L$ is not a sporadic simple group by~\cite{Giudici2006}.

We thus have $L=\A_n$ acting naturally on a set $\Omega$ of $n$ points. Since $G$ has a subgroup $H$ which has at least two nonsolvable composition factors, we know that $n\geqslant10$. Then according to Theorem~D and its Remark~2 in~\cite{LPS1990}, one of the following cases appears.
\begin{itemize}
\item[(i)] $H$ is $k$-homogeneous on $\Omega$ and $\A_{n-k}\leqslant K\leqslant\Sy_{n-k}\times\Sy_k$ for some $1\leqslant k\leqslant5$.
\item[(ii)] $\A_{n-k}\leqslant H\leqslant\Sy_{n-k}\times\Sy_k$ and $K$ is $k$-homogeneous on $\Omega$ for some $1\leqslant k\leqslant5$.
\item[(iii)] $n=10$, $H$ is a transitive subgroup of $\Sy_{10}$ such that $\A_5\times\A_5\vartriangleleft H\leqslant\Sy_5\wr\Sy_2$, and $K=\SL_2(8)$ or $\SL_2(8).3$.
\end{itemize}
Note that all the $k$-homogeneous permutation groups with $k\geqslant2$ are known: the $k$-transitive permutation groups are listed, for example in~\cite[Tables~7.3 and~7.4]{Cameron1999}, and the $k$-homogeneous but not $k$-transitive groups are classified in~\cite{Kantor1972}.

First assume case~(i) appears. If $k\geqslant2$, then the classification of $k$-homogeneous permutation groups shows that $H$ has at most one nonsolvable composition factor, a contradiction. Hence $k=1$ as in part~(a).

Next assume that case~(ii) appears. Then $k=5$ because $H$ has at least two nonsolvable composition factors. As $K\ngeqslant\A_n$, we conclude from the classification of $5$-homogeneous permutation groups~\cite{Cameron1999,Kantor1972} that $(n,K)=(12,\M_{12})$ or $(24,\M_{24})$. This is described in part~(c) or~(d), respectively, of Theorem~\ref{thm1}.

Finally, case~(iii) leads to part~(b). Thus the lemma is true.
\end{proof}

To complete the proof of Theorem~\ref{thm1}, we may assume by Lemma~\ref{lem2} that $L$ is a classical group of Lie type and assume by Lemmas~\ref{Maximality}--\ref{Embedding} that maximal subgroups of $G$ containing $H$ and $K$ respectively are core-free in $G$. We summarize these assumptions in the following hypothesis for later convenience.

\begin{hypothesis}\label{hyp1}
Let $G$ be an almost simple group with socle $L$ classical of Lie type. Suppose that $G=HK$ is a factorization of $G$ with $H$ having at least two nonsolvable composition factors and $K$ core-free, and $A$ and $B$ are core-free maximal subgroups of $G$ containing $H$ and $K$, respectively.
\end{hypothesis}

Notice that $A$ has at least one nonsolvable composition factor since the subgroup $H$ of $A$ is nonsolvable. In Section~5, we deal with the case where $A$ has at least two nonsolvable composition factors. Then in Section~6, we deal with the case where $A$ has exactly one nonsolvable composition factor.

\section{At least two nonsolvable composition factors of $A$}

In this section we prove that under Hypothesis~\ref{hyp1}, if $A$ has at least two nonsolvable composition factors then $(L,H\cap L,K\cap L)$ lies in Table~\ref{tab1} or Table~\ref{tab5}, as part~(b) of Theorem~\ref{thm1} asserts.

For $\Sp_4(2^f)\leqslant G\leqslant\GaSp_4(2^f)$ with $f\geqslant2$, we have factorizations $G=XY$ with $X\cap\Soc(G)=\GO_4^-(2^f)$ or $\Sz(2^f)$ and $Y\cap L=\Sp_2(2^f)\wr\Sy_2$ (see~\cite[Theorem~A]{LPS1990}). The lemma below shows that, however, $G=XY$ does not hold if $Y$ is replaced by its index $2$ subgroup $\N_2[G]$.

\begin{lemma}\label{lem4}
Let $L=\Sp_4(2^f)$ with $f\geqslant2$, and $L\leqslant G\leqslant\GaSp_4(2^f)$. Suppose $G=X\N_2[G]$ for some subgroup $X$ of $G$. Then $X\geqslant L$.
\end{lemma}

\begin{proof}
Let $Y$ be a maximal subgroup of $G$ containing $T:=\N_2[G]$. We have $Y\cap L=\Sp_2(2^f)\wr\Sy_2$ and $G=XY$. It then follows from~\cite[Theorem~A]{LPS1990} that one of the following cases appears.
\begin{itemize}
\item[(i)] $X\geqslant L$.
\item[(ii)] $X\cap L\leqslant\GO_4^-(2^f)$.
\item[(iii)] $f$ is odd and $X\cap L\leqslant\Sz(2^f)$.
\end{itemize}
If case~(ii) appears, then the factorization $G=XT$ would imply that $\GaO_4^-(2^f)$ is transitive on the nonsingular $2$-dimensional symplectic subspaces, which is not true.

Suppose that case~(iii) appears. Let $Z$ be a maximal subgroup of $G$ containing $X$ such that $Z\cap L=\Sz(2^f)$. According to~\cite{Suzuki1962}, $Z\cap L$ contains an element $g$ whose matrix with respect to a standard basis $e_1,f_1,e_2,f_2$ for $\Sp_4(2^f)$ is
$$
\begin{pmatrix}
0&1&&\\
1&0&&\\
&&0&1\\
&&1&0
\end{pmatrix}.
$$
By the equivalence of statements~(a) and~(b) of Lemma~\ref{lem6}, we may assume that $T$ fixes $\langle e_1,f_1\rangle$. Thus $Z\cap T\geqslant\langle g\rangle$ and so $|Z\cap T|_2\geqslant2$. Since $G/L\leqslant f$ has odd order, we have $|Z|_2=|Z\cap L|_2$ and $|T|_2=|T\cap L|_2$. Moreover, the factorization $G=XT$ indicates $G=ZT$ since $X\leqslant Z$. Thus we derive that
\begin{eqnarray*}
|Z\cap T|_2&=&\frac{|Z|_2|T|_2}{|G|_2}=\frac{|Z\cap L|_2|T\cap L|_2}{|G|_2}\\
&=&\frac{|\Sz(2^f)|_2|\Sp_2(2^f)\times\Sp_2(2^f)|_2}{|\Sp_4(2^f)|_2}=1,
\end{eqnarray*}
a contradiction. Hence $X\geqslant L$ as the lemma asserts.
\end{proof}

We also need the following two technical lemmas.

\begin{lemma}\label{lem5}
Let $L=\Sp_{2\ell}(2^f)$ with $f\ell\geqslant2$, and $L\leqslant G\leqslant\GaSp_{2\ell}(2^f)$. Suppose $G=RX=RY$ for some nonsolvable subgroup $R$ and subgroups $X$ and $Y$ such that $X\cap L=\GO_{2\ell}^+(2^f)$ and $Y\cap L=\GO_{2\ell}^-(2^f)$. Then either $\Sp_{2\ell/k}(2^{fk})\leqslant R\cap L\leqslant\Sp_{2\ell/k}(2^{fk}).k$ for some integer $k$ dividing $\ell$, or $\G_2(2^{f\ell/3})\leqslant R\cap L\leqslant\G_2(2^{f\ell/3}).(\ell/3)$ with $\ell$ divisible by $3$.
\end{lemma}

\begin{proof}
We use induction on $\ell$. First suppose $\ell=1$. Take any $r\in\ppd(2,2f)$. Then $r$ divides $|G|/|X|$, and so $r$ divides $|R|$ since $G=RX$. As $r>f$ and $|R|/|R\cap L|\leqslant f$, it follows that $r$ divides $|R\cap L|$. Thus, from the classification of subgroups of $\Sp_2(2^f)=\SL_2(2^f)$ (see for example~\cite[II~\S8]{Huppert1967}) we conclude that $R\cap L=\Sp_2(2^f)$ as $R\cap L$ is an nonsolvable subgroup of $\Sp_2(2^f)$ of order divisible by $r$.

Next suppose that $\ell>1$ and the lemma holds for each smaller $\ell$. Let $T$ be a maximal subgroup of $G$ containing $R$, and assume that $T$ is core-free in $G$ in view of Lemma~\ref{Embedding}. It follows that $G=TX=TY$. Then by~\cite[Theorem~A]{LPS1990}, either $T\cap L=\Sp_{2a}(2^{fb}).b$ with $ab=\ell$ and $b$ prime, or $T\cap L=\G_2(2^f)$ with $\ell=3$.

\underline{Case 1.} $T\cap L=\Sp_{2a}(2^{fb}).b$, where $ab=\ell$ and $b$ is a prime divisor of $\ell$. In this case, $T$ is an almost simple group with socle $\Sp_{2a}(2^{fb})$, and we deduce from Lemma~\ref{lem9} that $T=R(X\cap T)=R(Y\cap T)$. Furthermore, it is clear that $X\cap\Soc(T)=\GO_{2a}^+(2^{fb})$ and $Y\cap\Soc(T)=\GO_{2a}^-(2^{fb})$. Then by our inductive hypothesis, either
\begin{itemize}
\item[(i)] $\Sp_{2a/d}(2^{fbd})\leqslant R\cap\Soc(T)\leqslant\Sp_{2a/d}(2^{fbd}).d$ for some divisor $d$ of $a$, or
\item[(ii)] $\G_2(2^{fba/3})\leqslant R\cap\Soc(T)\leqslant\G_2(2^{fba/3}).(a/3)$ with $a$ divisible by $3$.
\end{itemize}
Note $R\cap\Soc(T)\leqslant R\cap L\leqslant(R\cap\Soc(T)).((T\cap L)/\Soc(T))=(R\cap\Soc(T)).b$. For~(i), we obtain by writing $k=bd$ that
$$
\Sp_{\ell/k}(2^{fk})\leqslant R\cap L\leqslant(\Sp_{\ell/k}(2^{fk}).d).b=\Sp_{\ell/k}(2^{fk}).k.
$$
For~(ii), we have
$$
\G_2(2^{f\ell/3})\leqslant R\cap L\leqslant(\G_2(2^{f\ell/3}).(a/3)).b=\G_2(2^{f\ell/3}).(\ell/3).
$$
Hence the conclusion of the lemma holds.

\underline{Case 2.} $\ell=3$ and $T\cap L=\G_2(2^f)$. In this case, $\G_2(2^f)\leqslant T\leqslant\G_2(2^f).f$, and $T=R(X\cap T)$ with $X\cap T\cap L=\SL_3(2^f).2$ (see~\cite[5.2.3(b)]{LPS1990}). Since $R$ is nonsolvable, we conclude $R\geqslant\Soc(T)$ by~\cite{HLS1987}, which implies $R\cap L=\G_2(2^f)$. This satisfies the conclusion of the lemma. The proof is thus completed.
\end{proof}

\begin{lemma}\label{lem10}
Let $G$ be an almost simple group with socle $L=\PSp_{2\ell}(q)$, where $\ell\geqslant3$ and $q\geqslant4$. Suppose $G=X\N_2[G]$ for some subgroup $X$ of $G$. Then one of the following holds.
\begin{itemize}
\item[(a)] $X\geqslant L$.
\item[(b)] $\ell=3$, $q$ is even and $X\cap L=\G_2(q)$.
\item[(c)] $\ell$ is even, $q=4$ and $\Sp_\ell(q^2)\leqslant X\cap L\leqslant\Sp_\ell(q^2).2$.
\item[(d)] $\ell=6$, $q=4$ and $\G_2(q^2)\leqslant X\cap L\leqslant\G_2(q^2).2$.
\end{itemize}
\end{lemma}

\begin{proof}
If $X$ is not core-free in $G$, then part~(a) of the lemma holds. In the remaining of the proof we assume that $X$ is core-free in $G$. Let $T$ be a maximal subgroup of $G$ containing $X$ and let $Y=\N_2[G]$. It follows from $G=XY$ that $G=TY$. By Lemma~\ref{Embedding} we may assume that $T$ is core-free in $G$. Then appealing to~\cite[Theorem~A]{LPS1990} we have the following two cases.
\begin{itemize}
\item[(i)] $\ell=3$, $q$ is even and $T\cap L=\G_2(q)$.
\item[(ii)] $\ell$ is even, $q=4$ and $T\cap L=\Sp_\ell(q^2).2$.
\end{itemize}
Note that $T=X(T\cap Y)$ by Lemma~\ref{lem9}.

First assume that case~(i) appears. Here $T$ is almost simple with socle $\G_2(q)$ and $T\cap Y\cap L=\Sp_2(q)\times\Sp_2(q)$ (see~\cite[5.2.3(b)]{LPS1990}). According to~\cite{HLS1987}, the factorization $T=X(T\cap Y)$ implies that $X\geqslant\Soc(T)$. Consequently, $X\cap L=\G_2(q)$ as described in part~(b) of the lemma.

Next assume that case~(ii) appears. In this case, $\Sp_\ell(q^2)\leqslant T\leqslant\GaSp_\ell(q^2)$ and $T\cap Y\leqslant\N_2[T]$ (see~\cite[3.2.1(a)]{LPS1990}). If $X$ is not core-free in $T$, then part~(c) of the lemma holds. Now assume that $X$ is core-free in $T$. Then by~\cite[Theorem~A]{LPS1990}, $\ell=6$ and $X\cap\Soc(T)\leqslant\G_2(q)$. Applying the same argument as in case~(i) (with $G$ therein replaced by $T$ here) to the factorization $T=X\N_2[T]$ we conclude that $X\cap\Soc(T)=\G_2(q)$. This implies $\G_2(q^2)\leqslant X\cap L\leqslant\G_2(q^2).2$, as in part~(d) of the lemma.
\end{proof}

Now we can prove the main result of this section.

\begin{lemma}\label{prop1}
Suppose \emph{Hypothesis~\ref{hyp1}} and that $A$ has at least two nonsolvable composition factors, then $(L,H\cap L,K\cap L)$ lies in \emph{Table~\ref{tab1}} or \emph{Table~\ref{tab5}}.
\end{lemma}

\begin{proof}
By Lemma~\ref{lem1}, $(L,A\cap L,B\cap L)=(L,X,Y)$ as in Table~\ref{tab2}. Let $G=L.\calO$ with $\calO\leqslant\Out(L)$. We discuss the eight rows of Table~\ref{tab2} in order.

\underline{Row~1.} In this row, $\calO\leqslant f$ and $A\cap L=(C\times D).2$ with $C\cong D=\Sp_{2\ell}(2^f)$. Let $M=(C\times D).\calO=\N_{2\ell}[G]$ be the subgroup of index $2$ in $A$ such that $M\cap L=C\times D$. Since $A\cap B=(\GO_{2\ell}^+(2^f)\times\GO_{2\ell}^-(2^f)).\calO$ (see~\cite[3.2.4(b)]{LPS1990}), we have $A\cap B\leqslant M$. Applying Lemma~\ref{lem9}, we deduce from $G=HB$ that $A=H(A\cap B)$ and then $M=(H\cap M)(A\cap B)$. Note that $M\leqslant(C.\calO)\times(D.\calO)$. Let $\pi_+$ and $\pi_-$ be the projections from $M$ to $C.\calO$ and $D.\calO$ respectively. Then $M^{\pi_\varepsilon}=(H\cap M)^{\pi_\varepsilon}(A\cap B)^{\pi_\varepsilon}$ for $\varepsilon=\pm$.

If $H\leqslant M$, then $G=\N_{2\ell}[G]B$, which implies that $B=\GO_{4\ell}^-(2^f).\calO$ is transitive on the nonsingular symplectic subspaces of dimension $2\ell$, a contradiction. Hence there exists $g\in H$ such that $C^g=D$ and so $(H\cap C)^g=H\cap D$. As a consequence, $(H\cap(C\times D))^{\pi_+}\cong(H\cap(C\times D))^{\pi_-}=U$ for some $U\leqslant\Sp_{2\ell}(2^f)$, and $(H\cap M)^{\pi_\varepsilon}\leqslant U.\calO$ for $\varepsilon=\pm$. Since $M^{\pi_\varepsilon}=(H\cap M)^{\pi_\varepsilon}(A\cap B)^{\pi_\varepsilon}$ and $(A\cap B)^{\pi_\varepsilon}=\GO_{2\ell}^\varepsilon(2^f).\calO$ for $\varepsilon=\pm$, we conclude from Lemma~\ref{lem5} that either $\Sp_{2\ell/k}(2^{fk})\leqslant U\leqslant\Sp_{2\ell/k}(2^{fk}).k$ for some integer $k$ dividing $\ell$, or $\G_2(2^{f\ell})\leqslant U\leqslant\G_2(2^{f\ell/3}).(\ell/3)$ with $\ell$ divisible by $3$. In particular, $U$ is an almost simple group. This implies $H\cap(C\times D)\geqslant\Soc(U)\times\Soc(U)$ as $H$ has at least two nonsolvable composition factors. It follows that
$$
H\cap L=(H\cap(C\times D)).2=(\Soc(U)\times\Soc(U)).R.2,
$$
where either $\Soc(U)=\Sp_{2\ell/k}(2^{fk})$ and $R\leqslant k\times k$ with $k$ dividing $\ell$, or $\Soc(U)=\G_2(2^{f\ell/3})$ and $R\leqslant(\ell/3)\times(\ell/3)$ with $\ell$ divisible by $3$. Next we show that $K\cap L\geqslant\Omega_{4\ell}^-(2^f)$ and thus either row~1 or row~2 of Table~\ref{tab1} appears.

Suppose $K\cap L\ngeqslant\Omega_{4\ell}^-(2^f)$. Then we have a factorization $B=(A\cap B)K$ with $K$ core-free in $B$. Moreover, $B$ is an almost simple group with socle $\Omega_{2\ell}^-(2^f)$, and $A\cap B=(\GO_{2\ell}^+(2^f)\times\GO_{2\ell}^-(2^f)).\calO$ is a maximal subgroup of $B$. However, inspection of~\cite[Tables~1--3]{LPS1990} shows that no such factorization exists. This contradiction proves that $K\cap L\geqslant\Omega_{4\ell}^-(2^f)$, as desired.

\underline{Row~2.} From~\cite[Table~1]{LPS1990} we see that $G=L.2$ and $A=(C\times D).2<(C.2)\times(D.2)$, where $C=\Sp_2(4)$ and $D=\Sp_{4\ell-2}(4)$. Let $\varphi$ be the projection from $A$ to $D.2$. Then since $A=H(A\cap B)$ by Lemma~\ref{lem9}, we have $A^\varphi=H^\varphi(A\cap B)^\varphi$. Moreover, $A^\varphi$ is an almost simple group with socle $\Sp_{4\ell-2}(4)$, and $(A\cap B)^\varphi\leqslant\N_2[A^\varphi]$. Suppose that $H^\varphi$ is core-free in $A^\varphi$. Then as $A^\varphi=H^\varphi\N_2[A^\varphi]$, we conclude from~\cite[Theorem~A]{LPS1990} that $\ell=2$ and $H^\varphi\leqslant\G_2(4).2$. Consequently, $H\cap L\leqslant C\times M$ with $M=\G_2(4)<D$. Applying Lemma~\ref{lem8} to the factorization $G=HB$ we see that $|H\cap L|$ is divisible by $|L|/|B|=|\Sp_8(4)|/|\Sp_4(16).4|=4|\G_2(4)|=4|M|$. Since the smallest index of proper subgroups of $M=\G_2(4)$ is $416$ (see~\cite{atlas}), which is larger than $|C|/4$, it follows that $H\cap L\geqslant M$. This in turn forces $H\cap L=C\times M$ in order that $H$ has at least two nonsolvable composition factors. However, $A\cap B\geqslant C$ (see~\cite[3.2.1(a)]{LPS1990}), whence
$$
|H\cap B|\geqslant|C|>\frac{|C|}{4}=\frac{|C||M||B|}{|L|}=\frac{|H\cap L||B|}{|L|}\geqslant\frac{|H||B|}{|G|},
$$
contradicting the factorization $G=HB$. Thus, $H^\varphi$ is not core-free in $A^\varphi$, and so $H\cap L\geqslant D$. As $H$ has at least two nonsolvable composition factors, we thereby obtain
$$
H\cap L=C\times D=\Sp_2(4)\times\Sp_{4\ell-2}(4).
$$
If $K\cap L\geqslant\Sp_{2\ell}(16)$, then row~3 of Table~\ref{tab1} appears. Now assume $K\cap L\ngeqslant\Sp_{2\ell}(16)$. Note that $B$ is almost simple with socle $\Sp_{2\ell}(16)$ and $A\cap B=\Sp_2(4)\times\Sp_{2\ell-2}(16)<\N_2[B]$. We have $B=\N_2[B]K$ with $K$ core-free in $B$. Thereby we conclude from Lemma~\ref{lem10} that $\ell=3$ and $\G_2(16)\leqslant K\cap L\leqslant\G_2(16).2$. This leads to row~8 of Table~\ref{tab1}.

\underline{Row~3.} In this row, $\calO\leqslant f$ and $A\cap L=(C\times D).2$ with $C\cong D=\Sp_2(2^f)$. Applying Lemma~\ref{lem8} to the factorization $G=HB$ we see that $f|H\cap L|$ is divisible by $|L|/|B\cap L|=2^{2f}(2^{2f}-1)(2^f+1)$. Consequently, as $|H\cap L|/|H\cap(C\times D)|$ divides $|A\cap L|/|C\times D|=2$, $f|H\cap(C\times D)|$ is divisible by $2^{2f-1}(2^{2f}-1)(2^f+1)$. Moreover, $|H\cap(C\times D)|$ divides $|H\cap C||D|=|H\cap C||\Sp_2(2^f)|$. It follows that $f|H\cap C|$ is divisible by $2^{2f-1}(2^{2f}-1)(2^f+1)/|\Sp_2(2^f)|=2^{f-1}(2^f+1)$. Then by the classification of subgroups of $\Sp_2(2^f)$ we conclude that $H\cap C=\Sp_2(2^f)$. For the same reason, $H\cap D=\Sp_2(2^f)$. Therefore, $H\cap L\geqslant C\times D$, which means that
$$
H\cap L=(C\times D).P=(\Sp_2(2^f)\times\Sp_2(2^f)).P
$$
with $P\leqslant2$. Since $G=AK$, by Lemma~\ref{lem9} we have $B=(A\cap B)K$. Note that $B$ is almost simple with socle $\Sz(2^f)$ and $A\cap B=\D_{2(2^f-1)}.\calO$ (see~\cite[5.1.7(b)]{LPS1990}). We then infer from~\cite{HLS1987} that $K\geqslant\Soc(B)$. Hence $K\cap L=\Sz(2^f)$, as in row~6 of Table~\ref{tab1}.

\underline{Row~4.} Let $N$, $M$, $\phi$, $A$, $\pi_1$ and $\pi_2$ be as in the proof of Example~\ref{exa4}. Since $G$ has a unique conjugacy class of maximal subgroups isomorphic to $B$, we may assume that $B=\langle N,\phi^e\rangle$ by Lemma~\ref{lem6}(b), where $e=|\calO|$. Then $A\cap B=M\rtimes\langle\phi^e\rangle$. By Lemma~\ref{lem9}, we derive from $G=HB$ and $H\leqslant A$ that $A=H(A\cap B)$, so $A^{\pi_2}=H^{\pi_2}(A\cap B)^{\pi_2}$. Note that $A^{\pi_2}=\Sp_4(2^f)$ has a graph automorphism fusing its $\calC_2$-subgroups with $\calC_8$ subgroups, and $M^{\pi_2}=\Omega_4^+(2^f)$. We thereby obtain a factorization $A^{\pi_2\gamma}=H^{\pi_2\gamma}\N_2[A^{\pi_2\gamma}]$ for some automorphism $\gamma$ of $A^{\pi_2}$. Hence by Lemma~\ref{lem4}, $H^{\pi_2}\geqslant\Soc(A^{\pi_2})=\Sp_4(2^f)$, and so
$$
H\cap L=U\times\Sp_4(2^f)
$$
for some subgroup $U$ of $(A\cap L)^{\pi_1}=\Sp_2(2^f)$. Moreover, $U$ is nonsolvable because $H\cap L$ has at least two nonsolvable composition factors. Next we show that $K\cap L=\G_2(2^f)$ and thus $(L,H\cap L,K\cap L)$ lies in row~7 of Table~\ref{tab1}.

Suppose $K\cap L<B\cap L=\G_2(2^f)$. Then is $K$ core-free in $B$. By Lemma~\ref{lem9}, we derive from $G=AK$ and $K\leqslant B$ that $B=(A\cap B)K$. However, as $B$ is an almost simple group with socle $\G_2(2^f)$ and $A\cap B=(\Sp_2(2^f)\times\Sp_2(2^f)).\calO$, it is seen in~\cite{HLS1987} that no such factorization $B=(A\cap B)K$ exists. This contradiction yields $K\cap L=\G_2(2^f)$, as desired.

\underline{Rows~5 and~6.} In these two rows, $L=\POm_{4\ell}^+(q)$ with either $\ell\geqslant3$ and $q\geqslant4$ or $\ell=2$ and $q\geqslant5$. Putting $c=\gcd(2,\ell,q-1)$, we have $A\cap L=(C\times D).c$ with $C=\PSp_2(q)$ and $D=\PSp_{2\ell}(q)$. Let $M$ be the subgroup of index $c$ in $A$ such that $M\cap L=C\times D$. Then $A\cap B\leqslant M$ and $M\leqslant(C.\calO)\times(D.\calO)$. Denote by $\varphi$ the projection from $M$ to $D.\calO$. Applying Lemma~\ref{lem9}, we deduce from $G=HB$ that $A=H(A\cap B)$ and then $M=(H\cap M)(A\cap B)$. Thereby we have $M^\varphi=(H\cap M)^\varphi(A\cap B)^\varphi$. Moreover, $M^\varphi/\Cen_{M^\varphi}(D)$ is almost simple with socle $D=\PSp_{2\ell}(q)$ and $(A\cap B)^\varphi\leqslant\N_2[M^\varphi]$ (see~\cite[3.6.1(d)]{LPS1990}). It follows that $M^\varphi=(H\cap M)^\varphi\N_2[M^\varphi]$, and so one of the following four cases appears by Lemma~\ref{lem10} .
\begin{itemize}
\item[(i)] $(H\cap M)^\varphi\geqslant D$.
\item[(ii)] $\ell=3$, $q$ is even and $(H\cap M)^\varphi\cap D=\G_2(q)$.
\item[(iii)] $\ell$ is even, $q=4$ and $\Sp_\ell(q^2)\leqslant(H\cap M)^\varphi\cap D\leqslant\Sp_\ell(q^2).2$.
\item[(iv)] $\ell=6$, $q=4$ and $\G_2(q^2)\leqslant(H\cap M)^\varphi\cap D\leqslant\G_2(q^2).2$.
\end{itemize}
Note that $H\cap L$ has at least two nonsolvable composition factors as $H$ has at least two nonsolvable composition factors. If case~(i) appears, then $H\cap L=(U\times\PSp_{2\ell}(q)).P$ with $P\leqslant\gcd(2,\ell,q-1)$, where $U$ is an nonsolvable subgroup of $\PSp_2(q)$. If case~(ii) appears, then $H\cap L=U\times\G_2(q)$ with $\ell=3$ and $q$ even, where $U$ is an nonsolvable subgroup of $\Sp_2(q)$. For case~(iii), $H\cap L=(\Sp_2(4)\times\Sp_\ell(16)).P$ with $\ell$ even, $q=4$ and $P\leqslant2$. For case~(iv) we have $H\cap L=(\Sp_2(4)\times\G_2(16)).P$ with $\ell=6$, $q=4$ and $P\leqslant2$. Next we show that $K\cap L=\Omega_{4\ell-1}(q)$ and thus cases~(i)--(iv) lead to rows~1, 4, 3, 6, respectively, of Table~\ref{tab5}.

Suppose $K\cap L\neq B\cap L=\Omega_{4\ell-1}(q)$. From $G=AK$ we derive that $|K|$ is divisible by $|G|/|A|$, whence each prime in $\ppd(q,4\ell-2)\cup\ppd(q,2\ell-1)\cup\ppd(q,4\ell-4)$ divides $|K|$. However, by Lemma~\ref{lem9} we have $B=(A\cap B)K$, and \cite[Theorem~A]{LPS1990} shows that no such group $K$ satisfies these conditions. This contradiction implies $K\cap L=\Omega_{4\ell-1}(q)$, as desired.

\underline{Row~7.} Here $L=\Omega_8^+(2)$, $H\cap L\leqslant(\SL_2(4)\times\SL_2(4)).2^2$ and $K\cap L\leqslant\Sp_6(2)$. In order that $H\cap L$ has at least two nonsolvable composition factors, we have $H\cap L=(\SL_2(4)\times\SL_2(4)).P$ with $P\leqslant2^2$. Furthermore, computation in \magma~\cite{magma} shows $K\cap L=\Sp_6(2)$ in order that $G=HK$. Thus row~2 (with $\ell=1$) of Table~\ref{tab5} appears.

\underline{Row~8.} In this row, $L=\Omega_8^+(4)$, $A\cap L=(\SL_2(16)\times\SL_2(16)).2^2$ and $B\cap L\leqslant\Sp_6(4)$. Then $|\calO|$ divides $4$, and thus Lemma~\ref{lem8} implies that $|L|$ divides $4|K\cap L||H\cap L|$. In particular, $|L|$ divides $4|B\cap L||H\cap L|$, which indicates that $|H\cap L|$ is divisible by $17$. Then as $H\cap L$ is a subgroup of $A\cap L=(\SL_2(16)\times\SL_2(16)).2^2$ with at least two nonsolvable composition factors, we conclude that either
$$
H\cap L=(\SL_2(4)\times\SL_2(16)).P
$$
with $P\leqslant2$, or
$$
H\cap L=(\SL_2(16)\times\SL_2(16)).P
$$
with $P\leqslant2^2$. For the former, since $4|K\cap L|$ is divisible by $|L|/|(\SL_2(4)\times\SL_2(16)).2|=2^{17}\cdot3^3\cdot5^2\cdot7\cdot13\cdot17$, we have $K\cap L=\Sp_6(4)$, which leads to row~3 (with $c=1$ and $\ell=1$) of Table~\ref{tab5}. Now assume that $H\cap L=(\SL_2(16)\times\SL_2(16)).P$ with $P\leqslant2^2$. It follows that $4|K\cap L|$ is divisible by $|L|/|(\SL_2(16)\times\SL_2(16)).2^2|=2^{14}\cdot3^3\cdot5^2\cdot7\cdot13$. Therefore, $K\cap L=\Sp_6(4)$ or $\G_2(4)$. If $K\cap L=\G_2(4)$, then $|K|$ divides $|\Out(L)||K\cap L|=4|\G_2(4)|=|\Omega_8^+(4)|/|(\SL_2(16)\times\SL_2(16)).2^2|=|L|/|A\cap L|=|G|/|A|$, which implies that $A\cap K=1$. However, there exists no such factorization $G=AK$ according to~\cite[Theorem~1.1]{LPS2010}. Thus we have $K\cap L=\Sp_6(4)$, and so row~3 (with $c=2$ and $\ell=1$) of Table~\ref{tab5} appears.
\end{proof}

\section{Exactly one nonsolvable composition factor of $A$}

In this section we deal with the case where the group $A$ in Hypothesis~\ref{hyp1} has exactly one nonsolvable composition factor. For a group $X$, the \emph{solvable radical} of $X$, denoted by $\Rad(X)$, is the product of all the solvable normal subgroups of $X$. Note that if $X$ has exactly one nonsolvable composition factor then $X/\Rad(X)$ is almost simple. Hence from the factorization $A=H(A\cap B)$ we obtain a factorization
$$
A/\Rad(A)=(H\Rad(A)/\Rad(A))((A\cap B)\Rad(A)/\Rad(A))
$$
of the almost simple group $A/\Rad(A)$ with the factor $H\Rad(A)/\Rad(A)$ still having at least two nonsolvable composition factors. We shall see in Lemma~\ref{Intersection} that the other factor $(A\cap B)\Rad(A)/\Rad(A)$ is core-free.

For a group $X$, let $X^{(\infty)}$ be the smallest normal subgroup of $X$ such that $X/X^{(\infty)}$ is solvable.

\begin{lemma}\label{Intersection}
Suppose that $G$ is an almost simple group with socle classical of Lie type, and $G=AB$ with subgroups $A$ and $B$ maximal and core-free in $G$. If $A$ has exactly one nonsolvable composition factor and $R=\Rad(A)$, then $A/R$ is almost simple and $(A\cap B)R/R$ is core-free in $A/R$.
\end{lemma}

\begin{proof}
Let $S$ be the unique nonsolvable composition factor of $A$. Then $A/R$ is an almost simple group with socle $S$. Suppose for a contradiction that $(A\cap B)R/R$ contains $\Soc(A/R)=(A/R)^{(\infty)}=A^{(\infty)}R/R$. We deduce that $(A\cap B)R\geqslant A^{(\infty)}R$. Then since $|(A\cap B)R|$ divides $|A\cap B||R|$, it follows that $|A\cap B||R|$ is divisible by $|A^{(\infty)}R|$. From $G=AB$ we have
$$
\frac{|G|}{|B|}=\frac{|A|}{|A\cap B|}=\frac{|A||R|}{|A\cap B||R|}.
$$
Therefore, $|G|/|B|$ divides
$$
\frac{|A||R|}{|A^{(\infty)}R|}=\frac{|A||A^{(\infty)}\cap R|}{|A^{(\infty)}|}
=\frac{|A||\Rad(A^{(\infty)})|}{|A^{(\infty)}|}=\frac{|A|}{|A^{(\infty)}/\Rad(A^{(\infty)})|}.
$$
Since $A^{(\infty)}/\Rad(A^{(\infty)})$ is also an almost simple group with socle $S$, we see that $|A^{(\infty)}/\Rad(A^{(\infty)})|$ is divisible by $|S|$ and so $|G|/|B|$ divides $|A|/|S|$. However, one concludes from~\cite[Theorem~A]{LPS1990} that no such factorization $G=AB$ exists with $|G|/|B|$ dividing $|A|/|S|$. Hence $(A\cap B)R/R$ does not contain $\Soc(A/R)$, which means that $(A\cap B)R/R$ is core-free in $A/R$.
\end{proof}

The following lemma is needed later.

\begin{lemma}\label{lem20}
Let $\Omega_{4m}^+(4)\leqslant G\leqslant\GO_{4m}^+(4)$ with $m\geqslant2$, and $H$ be a subgroup of $G$ such that $G=H\N_1[G]$. If $H$ is contained in a $\calC_4$-subgroup of $G$, then $H\geqslant\Sp_{2m}(4)$.
\end{lemma}

\begin{proof}
Let $W$, $V$, $Q$ and $M$ be as in Lemma~\ref{lem17} with $q=4$, and regard $G$ as $\Omega(W\otimes V,Q)$ or $\GO(W\otimes V,Q)$ according to $G=\Omega_{4m}^+(4)$ or $\GO_{4m}^+(4)$, respectively. Let $A$ be the maximal $\calC_4$-subgroup of $G$ containing $H$. Then $G=A\N_1[G]$, and so by~\cite[Theorem~A]{LPS1990} we may regard $A$ as the subgroup of $G$ stabilizing $W\otimes V$. From~\cite{KL1990} we see that $A<\Soc(G)$ and thus $A=\Sp_2(4)\times M$. Denote the projection of $A$ onto $M$ by $\varphi$. Since $G=H\N_1[G]$ and $H\leqslant A$, we deduce from Lemma~\ref{lem9} that $A=H(A\cap\N_1[G])$, which yields $M=A^\varphi=H^\varphi(A\cap\N_1[G])^\varphi$. As $(A\cap\N_1[G])^\varphi\leqslant\N_2[M]$ (see~\cite[3.6.1(d)]{LPS1990}), we further deduce that $M=H^\varphi\N_2[M]$. According to~\cite[Theorem~A]{LPS1990}, $M=\Sp_{2m}(4)$ has no maximal factorization with a factor $\N_2[M]$. Hence $H^\varphi=M$. This implies $H\geqslant\Sp_{2m}(4)$ as $H\leqslant A=\Sp_2(4)\times\Sp_{2m}(4)$.
\end{proof}

Based on the results above, we are able to prove Theorem~\ref{thm1} by induction on the order of $G$. The inductive hypothesis is as follows.

\begin{hypothesis}\label{hyp2}
Suppose that for each almost simple group $G_1$ of order properly dividing $|G|$, if $G_1=H_1K_1$ is a factorization of $G_1$ with $H_1$ having at least two nonsolvable composition factors and $K_1$ core-free, then $(G_1,H_1,K_1)$ is described in Theorem~\ref{thm1}.
\end{hypothesis}

Suppose both Hypothesis~\ref{hyp1} and Hypothesis~\ref{hyp2} in the rest of this section. We first give the candidates for $(L,A\cap L,B\cap L)$ in the case where $A$ has exactly one nonsolvable composition factor.

\begin{lemma}\label{prop2}
Either $A$ has at least two nonsolvable composition factors, or the triple $(L,A\cap L,B\cap L)$ lies in \emph{Table~\ref{tab3}}.
\end{lemma}

\begin{table}[htbp]
\caption{Factorizations in Lemma~\ref{prop2} with $A$ having a unique nonsolvable composition factor}\label{tab3}
\centering
\begin{tabular}{|l|l|l|l|}
\hline
row & $L$ & $A\cap L$ & $B\cap L$\\
\hline
1 & $\Sp_{4\ell}(2^f)$,~$\ell\geqslant2$ & $\Sp_{4a}(2^{fb}).b$ with $ab=\ell$ and $b$ prime & $\GO_{4\ell}^-(2^f)$\\
2 & $\Sp_{8\ell}(2)$,~$\ell\geqslant1$ & $\GO_{8\ell}^+(2)$ & $\GO_{8\ell}^-(2)$\\
3 & $\Sp_{4\ell}(4)$,~$\ell\geqslant2$ & $\GO_{4\ell}^+(4)$ & $\GO_{4\ell}^-(4)$\\
\hline
4 & $\Omega_{8\ell}^+(2)$,~$\ell\geqslant2$ & $\Omega_{4\ell}^+(4).2^2$ & $\Sp_{8\ell-2}(2)$\\
5 & $\Omega_{8\ell}^+(4)$,~$\ell\geqslant2$ & $\Omega_{4\ell}^+(16).2^2$ & $\Sp_{8\ell-2}(4)$\\
6 & $\Omega_{8}^+(2^f)$,~$f\geqslant2$ & $\Sp_6(2^f)$ & $\Sp_6(2^f)$\\
\hline
\end{tabular}
\end{table}

\begin{proof}
The maximal factorizations $G=AB$ are classified by~\cite[Theorem~A]{LPS1990} and described in Tables~1--4 of~\cite{LPS1990}. Since $H$ is nonsolvable, $A$ is nonsolvable too. Thus, either $A$ has at least two nonsolvable composition factors, or $A$ has precisely one nonsolvable composition factor. Assume in the following that $A$ has precisely one nonsolvable composition factor.

Denote $R=\Rad(A)$, $\overline{A}=A/R$, $\overline{H}=HR/R$ and $\overline{A\cap B}=(A\cap B)R/R$. Then $\overline{A}$ is an almost simple group, and $\overline{H}$ has at least two nonsolvable composition factors. Moreover, applying Lemma~\ref{lem9} we derive from $G=HB$ and $H\leqslant A$ that $A=H(A\cap B)$. Hence $\overline{A}=\overline{H}\,\overline{A\cap B}$ with $\overline{A\cap B}$ core-free in $\overline{A}$ as Lemma~\ref{Intersection} asserts. Now by Hypothesis~\ref{hyp2}, either $\Soc(\overline{A})=\A_n$ for some $n\geqslant10$, or~$\Soc(\overline{A})$ lies in the second column of Table~\ref{tab1} or Table~\ref{tab5}. Checking this condition of $A$ for factorizations $G=AB$ in~Tables 1--4 of~\cite{LPS1990}, we deduce that either $(L,A\cap L,B\cap L)$ lies in Table~\ref{tab4} or one of the following cases appears.
\begin{itemize}
\item[(i)] $L=\Sp_{2m}(2^f)$, $A\cap L=\Sp_{2a}(2^{fb}).b$ with $ab=m$ and $b$ prime, and $B\cap L=\GO_{2m}^-(2^f)$.
\item[(ii)] $L=\Sp_{2m}(2)$, $A\cap L=\GO_{2m}^+(2)$, and $B\cap L=\GO_{2m}^-(2)$.
\item[(iii)] $L=\Sp_{2m}(4)$, $A\cap L=\GO_{2m}^+(4)$, and $B\cap L=\GO_{2m}^-(4)$.
\item[(iv)] $L=\Sp_6(2^f)$, $A\cap L=\Pa_1$, and $B\cap L=\G_2(2^f)$.
\item[(v)] $G=L=\Sp_8(2)$, $A=\Sy_{10}$, $B=\GO_8^-(2)$.
\item[(vi)] $L=\Omega_{2m}^+(2^f)$ with $m\geqslant6$ even, $A\cap L=\N_1$, and $B\cap L=\bighat\GU_m(q).2$.
\item[(vii)] $L=\Omega_{2m}^+(2^f)$ with $m\geqslant6$ even and $f\leqslant2$, $A\cap L=\Omega_m^+(4^f).2^2$, and $B\cap L=\N_1$.
\item[(viii)] $L=\Omega_{10}^-(2)$, $A\cap L=\A_{12}$, and $B\cap L=\Pa_1$.
\item[(ix)] $L=\Omega_{8}^+(2^f)$ with $f\geqslant2$, $A\cap L=\Omega_7(2^f)$, and $B\cap L=\Omega_7(2^f)$.
\item[(x)] $L=\POm_8^+(3)$, $A\cap L=\Omega_8^+(2)$, and $B\cap L=\Omega_7(3)$.
\end{itemize}

\begin{table}[htbp]
\caption{Factorizations other than cases~(i)--(x) in the proof of Lemma~\ref{prop2}}\label{tab4}
\centering
\begin{tabular}{|l|l|l|l|l|}
\hline
$L$ & $A\cap L$ & $B\cap L$ & condition & $e$\\
\hline
$\PSL_n(q)$ & $\Nor_L(\PSp_n(q))$ & $\Pa_1$ or~$\Pa_{n-1}$ & $n\geqslant4$ even and & $n$\\
 & & & $(n,q)\neq(6,2)$ & \\
$\PSL_n(q)$ & $\Nor_L(\PSp_n(q))$ & $\Stab(V_1\oplus V_{n-1})$ & $n\geqslant4$ even and & $n$\\
 & & & $(n,q)\neq(6,2)$ & \\
\hline
$\PSp_{2m}(q)$ & $\PSp_{2a}(q^b).b$ & $\Pa_1$ & $m\geqslant2a\geqslant4$ & $2m$\\
$\PSp_{2m}(q)$ & $\Pa_1$ & $\PSp_{2a}(q^b).b$ & $m\geqslant3$ and & $2m-2$\\
 & & & $(m,q)\neq(4,2)$ & \\
$\PSp_{2m}(2^f)$ & $\Sp_{2a}(2^{fb}).b$ & $\GO_{2m}^+(2^f)$ & $m\geqslant2a\geqslant4$ & $2m$\\
$\PSp_{2m}(2^f)$ & $\GO_{2m}^+(2^f)$ & $\Sp_{2a}(2^{fb}).b$ & $m\geqslant4$ & $2m-2$\\
$\PSp_{2m}(2^f)$ & $\Sp_m(4^f).2$ & $\N_2$ & $m\geqslant4$ even & $2m$\\
 & & & and $f\leqslant2$ & \\
$\PSp_{2m}(2)$ & $\N_2$ & $\Sp_m(4).2$ & $m\geqslant6$ even & $2m-2$\\
$\PSp_{2m}(4^f)$ & $\Sp_{2m}(2^f)$ & $\GO_{2m}^-(4^f)$ & $f\leqslant2$ and & $m$\\
 & & & $mf\neq3$ & \\
$\PSp_8(2)$ & $\GO_8^+(2)$ & $\PSL_2(17)$ & & $3$\\
\hline
$\PSU_{2m}(q)$ & $\Nor_L(\PSp_{2m}(q))$ & $\N_1$ & $(m,q)\neq(3,2)$ & $2m$\\
\hline
$\POm_{2m}^-(2^f)$ & $\N_1$ & $\bighat\GU_m(2^f)$ & $m\geqslant5$ odd & $2m-2$\\
$\POm_{2m}^-(2)$ & $\N_1$ & $\Omega_m^-(4).2$ & $m\geqslant6$ even & $2m-2$\\
$\POm_{2m}^-(4)$ & $\N_1$ & $\Omega_m^-(16).2$ & $m\geqslant4$ even & $2m-2$\\
\hline
$\POm_{2m}^+(2^f)$ & $\N_1$ & $\Pa_m$ or~$\Pa_{m-1}$ & $m\geqslant5$ & $2m-2$\\
$\POm_{2m}^+(2^f)$ & $\N_1$ & $\Nor_L(\Sp_2(2^f)\otimes\Sp_m(2^f))$ & $m\geqslant6$ even & $2m-2$\\
$\POm_{2m}^+(q)$ & $\Pa_1$ & $\bighat\GU_m(q).2$ & $m\geqslant6$ even & $2m-4$\\
$\POm_{2m}^+(2^f)$ & $\N_1$ & $\bighat\GL_m(q).2$ & $m\geqslant5$ & $2m-2$\\
$\POm_{2m}^+(2^f)$ & $\N_1$ & $\Omega_m^+(4^f).2^2$ & $m\geqslant6$ even & $2m-2$\\
 & & & and $f\leqslant2$ & \\
$\POm_{2m}^+(4)$ & $\N_2^+$ & $\bighat\GU_m(4).2$ & $m\geqslant6$ even & $2m-4$\\
$\POm_{16}^+(2^f)$ & $\Omega_9(2^f).a$ & $\N_1$ & $a\leqslant2$ & $8$\\
$\POm_{16}^+(2^f)$ & $\N_1$ & $\Omega_9(2^f).a$ & $a\leqslant2$ & $14$\\
$\POm_{24}^+(2)$ & $\N_1$ & $\Co_1$ & & $22$\\
\hline
$\POm_8^+(2^f)$ & $\Omega_7(2^f)$ & $\Pa_1$,~$\Pa_3$ or~$\Pa_4$ & $f\geqslant2$ & $6$\\
$\POm_8^+(2^f)$ & $\Omega_7(2^f)$ & $\bighat((2^f+1)\times\Omega_6^-(2^f)).2$ & $f\geqslant2$ & $3$\\
$\POm_8^+(2^f)$ & $\Omega_7(2^f)$ & $\bighat((2^f-1)\times\Omega_6^+(2^f)).2$ & $f\geqslant2$ & $6$\\
$\POm_8^+(4^f)$ & $\Omega_7(4^f)$ & $\Omega_8^-(2^f)$ & & $6$\\
$\POm_8^+(3)$ & $\Omega_8^+(2)$ & $\Pa_1$,~$\Pa_3$ or~$\Pa_4$ & & $6$\\
$\POm_8^+(3)$ & $\Omega_8^+(2)$ & $\Pa_{13}$,~$\Pa_{14}$ or~$\Pa_{34}$ & & $6$\\
$\Omega_8^+(4)$ & $\Omega_7(4)$ & $(\SL_2(16)\times\SL_2(16)).2^2$ & & $6$\\
\hline
\end{tabular}
\end{table}

Suppose that $(L,A\cap L,B\cap L)$ lies in Table~\ref{tab4}, so that $L=\PSL_n(q)$, $\PSp_{2m}(q)$, $\PSU_{2m}(q)$ or $\POm_{2m}^\pm(q)$. For each row in Table~\ref{tab4}, take $e$ to be the number as in the last column and take $r\in\ppd(q,e)$. Then one can verify that $r$ divides $|G|/|B|$ but not $|R|$. It follows that $r$ divides $H$ due to the factorization $G=HB$, and then $r$ divides $\overline{H}$. Further, $r$ divides $\overline{H}\cap\Soc(\overline{A})$ as $r$ does not divide $|\Out(\Soc(\overline{A}))|$. However, Hypothesis~\ref{hyp2} requires $\overline{H}\cap\Soc(\overline{A})$ to be in the third column of Table~\ref{tab1} or Table~\ref{tab5}, which implies that $|\overline{H}\cap\Soc(\overline{A})|$ is not divisible by $r$. This contradiction shows that $(L,A\cap L,B\cap L)$ does not lie in Table~\ref{tab4}. Thus, we only need to consider cases~(i)--(x) in the following.

First assume that case~(i) appears. From Hypothesis~\ref{hyp2} we deduce that either $a$ is even, or $a=3$ and $\overline{H}\cap\Soc(\overline{A})\leqslant\Sp_4(2^{fb})\times\Sp_2(2^{fb})$. The former immediately leads to row~1 of Table~\ref{tab3}, so we assume the latter now. Take $r\in\ppd(2,3fb)$. Then as $r$ is coprime to both $|\overline{H}|$ and $|R|$, $r$ is coprime to $|H|$. However, $r$ divides $|\Sp_{6b}(2^f)|/|\GO_{6b}^-(2^f)|=|L|/|B\cap L|$, contradicting the factorization $G=HB$ by Lemma~\ref{lem8}.

For case~(ii), we deduce from Hypothesis~\ref{hyp2} that $m$ is divisible by $4$, and so row~2 of Table~\ref{tab3} appears.

For case~(iii), we deduce from Hypothesis~\ref{hyp2} that $m$ is even, and so row~3 of Table~\ref{tab3} appears.

Suppose that case~(iv) occurs. Then either $f\geqslant2$ and $\Soc(\overline{A})=\Sp_4(2^f)$, or $f=1$ and $\Soc(\overline{A})=\A_6$. From Hypothesis~\ref{hyp2} we deduce that $f\geqslant2$ and any prime $r$ in $\ppd(2,4f)$ is coprime to $|\overline{H}|$. It follows that $r$ does not divide $|H|$ since $r$ does not divide $|R|$. However, $r$ divides $|\Sp_6(2^f)|/|\G_2(2^f)|=|L|/|B\cap L|$, contradicting the factorization $G=HB$ by Lemma~\ref{lem8}.

Suppose that case~(v) occurs. Since $G=HB$ and $H\leqslant A$, we have the factorization $A=H(A\cap B)$ by Lemma~\ref{lem9}. However, $A\cap B=\Sy_7\times\Sy_3$ by~\cite[5.1.9]{LPS1990}, which implies that the factorization $A=H(A\cap B)$ does not satisfy Theorem~\ref{thm1}. This is contrary to Hypothesis~\ref{hyp2}.

Next suppose that case~(vi) occurs. Then as $2m-2\geqslant10$, we deduce from Hypothesis~\ref{hyp2} that $2m-2$ is divisible by $4$, contrary to the condition that $m$ is even.

For case~(vii), we deduce from Hypothesis~\ref{hyp2} that $m$ is divisible by $4$, whence row~4 or~5 of Table~\ref{tab3} appears.

Suppose that case~(viii) occurs. In this case, we have $(G,A)=(\Omega_{10}^-(2),\A_{12})$ or $(\GO_{10}^-(2),\Sy_{12})$. Since $G=HB$ and $H\leqslant A$, we have the factorization $A=H(A\cap B)$ by Lemma~\ref{lem9}. However, $A\cap B=\Sy_8\times\Sy_4$ by~\cite[5.2.16]{LPS1990}, which implies that the factorization $A=H(A\cap B)$ does not satisfy Theorem~\ref{thm1}. This is contrary to Hypothesis~\ref{hyp2}.

For case~(ix), in view of $\Omega_7(2^f)\cong\Sp_6(2^f)$ we see that row~6 of Table~\ref{tab3} appears.

Finally suppose that~(x) occurs. Since $G=HB$ and $H\leqslant A$, we have the factorization $A=H(A\cap B)$ by Lemma~\ref{lem9}. However, $A\cap B\cap L=2^6.\A_7$ by~\cite[5.1.15]{LPS1990}, which implies that the factorization $A=H(A\cap B)$ does not satisfy Theorem~\ref{thm1}. This is contrary to Hypothesis~\ref{hyp2}.
\end{proof}

In the forthcoming Lemmas~\ref{lem3}--\ref{lem15}, we analyze the six candidates for the triple $(L,A\cap L,B\cap L)$ in Table~\ref{tab3} to show that $(L,H\cap L,K\cap L)$ lies in Table~\ref{tab1} or Table~\ref{tab5}. Note that we have $A=H(A\cap B)$ and $B=(A\cap B)K$ from the factorization $G=HK$. Moreover, we see in Table~\ref{tab3} that both $A$ and $B$ have exactly one nonsolvable composition factor, so both $A/\Rad(A)$ and $B/\Rad(B)$ are almost simple.

\begin{lemma}\label{lem3}
Suppose $(L,A\cap L,B\cap L)$ lies in \emph{row~1} of \emph{Table~\ref{tab3}}. Then $(L,H\cap L,K\cap L)$ lies in \emph{Table~\ref{tab1}}.
\end{lemma}

\begin{proof}
Here $L=\Sp_{4\ell}(2^f)$, $A\cap L=\Sp_{4a}(2^{fb}).b$ and $B\cap L=\GO_{4\ell}^-(2^f)$ with $ab=\ell$ and $b$ prime. Since $A=H(A\cap B)$ and $A\cap B\cap L=\GO_{4a}^-(2^{fb}).b$, we infer from Hypothesis~\ref{hyp2} that $\Soc(A)$ is described in the second column of row~1, 2, 4, 5, 9 or~10 of Table~\ref{tab1}. Take any $r\in\ppd(2,4f\ell)$. Since $B=(A\cap B)K$ and $B$ is almost simple with socle $\Omega_{4\ell}^-(2^f)$, we derive from~\cite[Theorem~A]{LPS1990} that either $K\geqslant\Soc(B)$ or $K\leqslant\N_1[B]$. For the latter, $r$ divides $|H\cap L|$ because $r$ divides $|L|$ but not $|\N_1[B]|$.

\underline{Case~1.} Let $\Soc(A)=\Sp_{4a}(2^{fb})$ lie in the second column of row~1 of Table~\ref{tab1}. In this case, $H\cap\Soc(A)=(\Sp_{2a_0}(2^{fbb_0})\times\Sp_{2a_0}(2^{fbb_0})).R_0.2$ with $R_0\leqslant b_0\times b_0$, where $a_0$ and $b_0$ are positive integers such that $a=a_0b_0$. Let $a_1=a_0$ and $b_1=b_0b$. Then $a_1b_1=a_0b_0b=ab=\ell$ and $H\cap L=(\Sp_{2a_1}(2^{fb_1})\times\Sp_{2a_1}(2^{fb_1})).R.2$ with $R\leqslant b_1\times b_1$. In particular, $|H\cap L|$ is not divisible by $r$, and so $K\geqslant\Soc(B)$ by the conclusion of the previous paragraph. It follows that $K\cap L=\Omega_{4\ell}^-(2^f).Q$ with $Q\leqslant2$, and thus $(L,H\cap L,K\cap L)$ is as described in row~1 of Table~\ref{tab1}.

\underline{Case~2.} Let $\Soc(A)=\Sp_{4a}(2^{fb})$ lie in the second column of row~2 of Table~\ref{tab1}. Here $a=3a_1$ for some integer $a_1$ and $H\cap\Soc(A)=(\G_2(2^{fba_1})\times\G_2(2^{fba_1})).R_1.2$ with $R_1\leqslant a_1\times a_1$. Let $\ell_1=a_1b$. It follows that $\ell=ab=3a_1b=3\ell_1$, $L=\Sp_{12\ell_1}(2^f)$ and $H\cap L=(\G_2(2^{f\ell_1})\times\G_2(2^{f\ell_1})).R.2$ with $R\leqslant\ell_1\times\ell_1$. In particular, $|H\cap L|$ is not divisible by $r$, which implies that $K\geqslant\Soc(B)$. Thus, $K\cap L=\Omega_{4\ell}^-(2^f).Q=\Omega_{12\ell_1}^-(2^f).Q$ with $Q\leqslant2$, and so $(L,H\cap L,K\cap L)$ is as described in row~2 of Table~\ref{tab1}.

\underline{Case~3.} Let $\Soc(A)=\Sp_{4a}(2^{fb})$ lie in the second column of row~4 of Table~\ref{tab1}. Here $a\geqslant2$, $b=2$, $f=1$ and $H\cap\Soc(A)=(\Sp_2(4)\times\Sp_{2a}(4)).P$ with $P\leqslant2$. Thus, $L=\Sp_{4\ell}(2)=\Sp_{8a}(2)$. Since
$$
(H\cap L)/(H\cap\Soc(A))\cong(H\cap L)\Soc(A)/\Soc(A)\leqslant(A\cap L)/\Soc(A)=b,
$$
we have $\Sp_2(4)\times\Sp_{2a}(4)\leqslant H\cap L\leqslant(\Sp_2(4)\times\Sp_{2a}(4)).[4]$. In particular, $|H\cap L|$ is not divisible by $r$, and so $K\geqslant\Soc(B)$. Then we conclude that $K\cap L=\Omega_{4\ell}^-(2).Q=\Omega_{8a}^-(2).Q$ with $Q\leqslant2$ and $(L,H\cap L,K\cap L)$ is as described in row~5 of Table~\ref{tab1}.

\underline{Case~4.} Let$\Soc(A)=\Sp_{4a}(2^{fb})$ lie in the second column of row~5 or 10 of Table~\ref{tab1}. In this case we have $fb=1$, contrary to the condition that $b$ is prime.

\underline{Case~5.} Let $\Soc(A)=\Sp_{4a}(2^{fb})$ lie in the second column of row~9 of Table~\ref{tab1}. In this case, $a=3$, $b=2$, $f=1$ and $H\cap\Soc(A)=(\Sp_2(4)\times\G_2(4)).P$ with $P\leqslant2$. As $H\cap\Soc(A)\leqslant H\cap L\leqslant(H\cap\Soc(A)).b$, we derive that $H\cap L=(\Sp_2(4)\times\G_2(4)).P_1$, where $P_1\leqslant2^2$. It follows that $|H\cap L|$ is not divisible by $r$, and then $K\geqslant\Soc(B)$. Thus, $K\cap L=\Omega_{4\ell}^-(2).Q=\Omega_{24}^-(2).Q$ with $Q\leqslant2$, and so $(L,H\cap L,K\cap L)$ lies in row~10 of Table~\ref{tab1}.
\end{proof}

\begin{lemma}\label{lem11}
Suppose $(L,A\cap L,B\cap L)$ lies in \emph{row~2} of \emph{Table~\ref{tab3}}. Then $(L,H\cap L,K\cap L)$ lies in \emph{Table~\ref{tab1}}.
\end{lemma}

\begin{proof}
Here $G=\Sp_{8\ell}(2)$, $A=\GO_{8\ell}^+(2)$ and $B=\GO_{8\ell}^-(2)$. Since $A=H(A\cap B)$ and $A\cap B=\Sp_{8\ell-2}(2)\times2=\N_1[A]$ (see~\cite[3.2.4(e)]{LPS1990}), Hypothesis~\ref{hyp2} implies that $\Soc(A)$ is described in the second column of row~2 or~5 of Table~\ref{tab5}. Take any $r\in\ppd(2,8\ell-2)$. Since $B=(A\cap B)K$ and $B$ is almost simple with socle $\Omega_{8\ell}^-(2)$, we derive from~\cite[Theorem~A]{LPS1990} that either $K\geqslant\Soc(B)$ or $K\cap\Soc(B)\leqslant\Omega_{4\ell}^-(4).2$. For the latter, $r$ divides $|H\cap L|$ because $r$ divides $|L|$ but not $|\Omega_{4\ell}^-(4).2|$.

\underline{Case~1.} Let $\Soc(A)=\Omega_{8\ell}^+(2)$ lie in the second column of row~2 of Table~\ref{tab5}. In this case, $H\cap\Soc(A)=(\Sp_2(4)\times\Sp_{2\ell}(4)).P$ with $P\leqslant2^2$. Then since
$$
(H\cap L)/(H\cap\Soc(A))\cong(H\cap L)\Soc(A)/\Soc(A)\leqslant(A\cap L)/\Soc(A)=2,
$$
we have $\Sp_2(4)\times\Sp_{2\ell}(4)\leqslant H\cap L\leqslant(\Sp_2(4)\times\Sp_{2\ell}(4)).[8]$. In particular, $|H\cap L|$ is not divisible by $r$, which implies $K\geqslant\Soc(B)$ by the conclusion of the previous paragraph. Hence $K\cap L=\Omega_{8\ell}^-(2).Q$ with $Q\leqslant2$ and so $(L,H\cap L,K\cap L)$ lies in row~5 of Table~\ref{tab1}.

\underline{Case~2.} Let $\Soc(A)=\Omega_{8\ell}^+(2)$ lie in the second column of row~5 of Table~\ref{tab5}. In this case, $\ell=3$ and $H\cap\Soc(A)=(\Sp_2(4)\times\G_2(4)).P$ with $P\leqslant2^2$. Since $H\cap\Soc(A)\leqslant H\cap L\leqslant(H\cap\Soc(A)).2$, we have $\Sp_2(4)\times\G_2(4)\leqslant H\cap L\leqslant(\Sp_2(4)\times\G_2(4)).[8]$. In particular, $|H\cap L|$ is not divisible by $r$, and so $K\geqslant\Soc(B)$. It follows that $K\cap L=\Omega_{8\ell}^-(2).Q=\Omega_{24}^-(4).Q$ with $Q\leqslant2$, and thus $(L,H\cap L,K\cap L)$ lies in row~10 of Table~\ref{tab1}.
\end{proof}

\begin{lemma}\label{lem12}
Suppose $(L,A\cap L,B\cap L)$ lies in \emph{row~3} of \emph{Table~\ref{tab3}}. Then $(L,H\cap L,K\cap L)$ lies in \emph{Table~\ref{tab1}}.
\end{lemma}

\begin{proof}
Here $\ell\geqslant2$ and $(G,A,B)=(\GaSp_{4\ell}(4),\GaO_{4\ell}^+(4),\GaO_{4\ell}^-(4))$ according to~\cite[Theorem~A]{LPS1990}. Since $A=H(A\cap B)$ and $A\cap B=\N_1[A\cap L]=\Sp_{4\ell-2}(4)\times2$ (see~\cite[3.2.4(e)]{LPS1990}), we derive from Hypothesis~\ref{hyp2} that $\Soc(A)$ lies in the second column of row~1, 3, 4 or 6 of Table~\ref{tab5}. Take any $r\in\ppd(4,4\ell-2)$. Since $B=(A\cap B)K$ and $B$ is almost simple with socle $\Omega_{4\ell}^-(4)$, it follows from~\cite[Theorem~A]{LPS1990} that either $K\geqslant\Soc(B)$ or $K\cap\Soc(B)\leqslant\Omega_{2\ell}^-(16).2$. For the latter, $r$ divides $|H\cap L|$ because $r$ divides $|L|$ but not $|\Omega_{2\ell}^-(16).2|$.

\underline{Case~1.} Let $\Soc(A)=\Omega_{4\ell}^+(4)$ lie in the second column of row~1 of Table~\ref{tab5}. In this case, $H\cap\Soc(A)=U\times\Sp_{2\ell}(4)$, where $U$ is an nonsolvable subgroup of $\Sp_2(4)$. Therefore, $U=\Sp_2(4)$ and $H\cap\Soc(A)=\Sp_2(4)\times\Sp_{2\ell}(4)$. Since
$$
(H\cap L)/(H\cap\Soc(A))\cong(H\cap L)\Soc(A)/\Soc(A)\leqslant(A\cap L)/\Soc(A)=2,
$$
we deduce that $\Sp_2(4)\times\Sp_{2\ell}(4)\leqslant H\cap L\leqslant(\Sp_2(4)\times\Sp_{2\ell}(4)).2$. In particular, $|H\cap L|$ is not divisible by $r$, and so $K\geqslant\Soc(B)$ by the conclusion of the previous paragraph. Hence $K\cap L=\Omega_{4\ell}^-(4).Q$ with $Q\leqslant2$, and thus $(L,H\cap L,K\cap L)$ lies in row~4 of Table~\ref{tab1}.

\underline{Case~2.} Let $\Soc(A)=\Omega_{4\ell}^+(4)$ lie in the second column of row~3 or~6 of Table~\ref{tab5}. In this case we have $\ell=2\ell_1$ for some integer $\ell_1$. Moreover, either $H\cap\Soc(A)=(\Sp_2(4^c)\times\Sp_{2\ell_1}(16)).P$, or $\ell_1=3$ and $H\cap\Soc(A)=(\Sp_2(4^c)\times\G_2(16)).P$, where $c=1$ or $2$ and $P\leqslant2^2$. As $A\cap B=\N_1[A\cap L]<A\cap L$ (see~\cite[3.2.4(e)]{LPS1990}) and $A=H(A\cap B)$, we have a factorization
$$
A\cap L=(H\cap L)\N_1[A\cap L]
$$
of $A\cap L=\GO_{8\ell_1}^+(4)$. By Lemma~\ref{Embedding}, there exist groups $A^*$ and $H^*$ such that $\Omega_{8\ell_1}^+(4)\leqslant A^*\leqslant\GO_{8\ell_1}^+(4)$, $A^*=H^*\Soc(A)=H^*\N_1[A^*]$ and $H^*\cap\Soc(A)=H\cap\Soc(A)$. Let $M$ be a maximal subgroup of $A^*$ containing $H^*$. Then $A^*=M\N_1[A^*]$, and $M$ is core-free in $A^*$. From Hypothesis~\ref{hyp2} we see that $M$ cannot have a linear, unitary or plus type orthogonal group as its unique nonsolvable composition factor. Thus one of the following occurs by~\cite[Theorem~A]{LPS1990}.
\begin{itemize}
\item[(i)] $M\cap\Soc(A^*)=\Sp_2(4)\otimes\Sp_{4\ell_1}(4)$.
\item[(ii)] $M\cap\Soc(A^*)=\Omega_{4\ell_1}(16).2^2$.
\item[(iii)] $\ell_1=2$ and $M\cap\Soc(A^*)=\Sp_8(4).a$ with $a\leqslant2$.
\item[(iv)] $\ell_1=1$ and $M\cap\Soc(A^*)=\Sp_6(4)$.
\item[(v)] $\ell_1=1$ and $M\cap\Soc(A^*)=(\Sp_2(16)\times\Sp_2(16)).2^2$.
\end{itemize}
If (i) occurs, then $H^*\geqslant\Sp_{4\ell_1}(4)$ by Lemma~\ref{lem20}, which is not possible since we have either $H^*\cap\Soc(A)=(\Sp_2(4^c)\times\Sp_{2\ell_1}(16)).P$ or $\ell_1=3$ and $H^*\cap\Soc(A)=(\Sp_2(4^c)\times\G_2(16)).P$. If $M\cap\Soc(A^*)$ is as described in~(ii) or~(v), then the argument in~\cite[3.6.1(c)]{LPS1990} shows that $A^*\neq M\N_1[A^*]$, a contradiction. If~(iii) occurs, then we see from Hypothesis~\ref{hyp2} and the third column of Table~\ref{tab1} that $H^*$ cannot have $\Sp_2(4^c)\times\Sp_{2\ell_1}(16)=\Sp_2(4^c)\times\Sp_4(16)$ as a normal subgroup, a contradiction. Finally, suppose that~(iv) occurs. Then $M\cap\N_1[A^*]$ has a normal subgroup $\G_2(4)$ by~\cite[5.1.15]{LPS1990} while $H^*$ has a normal subgroup $\Sp_2(4^c)\times\Sp_2(16)$. By Lemma~\ref{lem9}, we deduce from $A^*=H^*\N_1[A^*]$ and $H^*\leqslant M$ that $M=H^*(M\cap\N_1[A^*])$. However, $M$ does not have such a factorization according to Hypothesis~\ref{hyp2}, a contradiction.

\underline{Case~3.} Let $\Soc(A)=\Omega_{4\ell}^+(4)$ lie in the second column of row~4 of Table~\ref{tab5}. In this case, $\ell=3$ and $H\cap\Soc(A)=(\Sp_2(4)\times\G_2(4))$. Since $H\cap\Soc(A)\leqslant H\cap L\leqslant(H\cap\Soc(A)).2$, we have $\Sp_2(4)\times\G_2(4)\leqslant H\cap L\leqslant(\Sp_2(4)\times\G_2(4)).2$. In particular, $|H\cap L|$ is not divisible by $r$, and so $K\geqslant\Soc(B)$. Hence $K\cap L=\Omega_{4\ell}^-(4).Q=\Omega_{12}^-(4).Q$ with $Q\leqslant2$, and then $(L,H\cap L,K\cap L)$ lies in row~9 of Table~\ref{tab1}.
\end{proof}

\begin{lemma}\label{lem13}
Suppose $(L,A\cap L,B\cap L)$ lies in \emph{row~4} of \emph{Table~\ref{tab3}}. Then $(L,H\cap L,K\cap L)$ lies in \emph{Table~\ref{tab5}}.
\end{lemma}

\begin{proof}
Here $L=\Omega_{8\ell}^+(2)$ with $\ell\geqslant2$, $A\cap L=\Omega_{4\ell}^+(4).2^2$ and $B\cap L=\Sp_{8\ell-2}(2)$. Furthermore, we have $G=L$ in order that $A$ is maximal in $G$. By Lemma~\ref{lem19}, $A\cap B=\Sp_{4\ell-2}(4)\times2<\GO_{4\ell}^+(4)$. Then as $A=H(A\cap B)$, we derive from Hypothesis~\ref{hyp2} that $\Soc(A)$ is described in the second column of row~1, 3, 4 or~6 of Table~\ref{tab5}. Since $B=(A\cap B)K$ and $|A\cap B|$ has no prime divisor in the set $\ppd(2,8\ell-2)\cup\ppd(2,8\ell-6)\cup\ppd(2,4\ell-1)$, we conclude that $|K|$ is divisible by each prime in $\ppd(2,8\ell-2)\cup\ppd(2,8\ell-6)\cup\ppd(2,4\ell-1)$. Hence by~\cite[Theorem~A]{LPS1990}, the factorization $B=(A\cap B)K$ implies that $K\geqslant\Soc(B)$, and so $K\cap L=\Sp_{8\ell-2}(2)$.

\underline{Case~1.} Let $\Soc(A)=\Omega_{4\ell}^+(4)$ lie in the second column of row~1 of Table~\ref{tab5}. In this case, $H\cap\Soc(A)=U\times\Sp_{2\ell}(4)$, where $U$ is an nonsolvable subgroup of $\Sp_2(4)$. Thus, $U=\Sp_2(4)$ and $H\cap\Soc(A)=\Sp_2(4)\times\Sp_{2\ell}(4)$. Since
$$
(H\cap L)/(H\cap\Soc(A))\cong(H\cap L)\Soc(A)/\Soc(A)\leqslant(A\cap L)/\Soc(A)=2^2,
$$
we deduce that $\Sp_2(4)\times\Sp_{2\ell}(4)\leqslant H\cap L\leqslant(\Sp_2(4)\times\Sp_{2\ell}(4)).2^2$. Therefore, $(L,H\cap L,K\cap L)$ lies in row~2 of Table~\ref{tab5}.

\underline{Case~2.} Let $\Soc(A)=\Omega_{4\ell}^+(4)$ lie in the second column of row~3 or~6 of Table~\ref{tab5}. Since $A=H(A\cap B)$ and $A\cap B=\Sp_{4\ell-2}(4)\times2\leqslant\N_1[\GO_{4\ell}^+(4)]$, we derive that
$$
\GO_{4\ell}^+(4)=(H\cap\GO_{4\ell}^+(4))(A\cap B)=(H\cap\GO_{4\ell}^+(4))\N_1[\GO_{4\ell}^+(4)].
$$
Then one obtains a contradiction along the same lines as in Case~2 of the proof of Lemma~\ref{lem12}.

\underline{Case~3.} Let $\Soc(A)=\Omega_{4\ell}^+(4)$ lie in the second column of row~4 of Table~\ref{tab5}. In this case, $\ell=3$ and $H\cap\Soc(A)=(\Sp_2(4)\times\G_2(4))$. As $H\cap\Soc(A)\leqslant H\cap L\leqslant(H\cap\Soc(A)).2^2$, we have $\Sp_2(4)\times\G_2(4)\leqslant H\cap L\leqslant(\Sp_2(4)\times\G_2(4)).2^2$. Thus, $(L,H\cap L,K\cap L)$ lies in row~5 of Table~\ref{tab5}.
\end{proof}

\begin{lemma}\label{lem14}
Suppose $(L,A\cap L,B\cap L)$ lies in \emph{row~5} of \emph{Table~\ref{tab3}}. Then $(L,H\cap L,K\cap L)$ lies in \emph{Table~\ref{tab5}}.
\end{lemma}

\begin{proof}
Here $L=\Omega_{8\ell}^+(4)$ with $\ell\geqslant2$, $A\cap L=\Omega_{4\ell}^+(16).2^2$ and $B\cap L=\Sp_{8\ell-2}(4)$. Since $A=H(A\cap B)$ and $A\cap B\cap L=\Sp_{4\ell-2}(16)$, we derive from Hypothesis~\ref{hyp2} that $\Soc(A)$ is described in the second column of row~1 or~4 of Table~\ref{tab5}. As $B=(A\cap B)K$ and $|A\cap B|$ has no prime divisor in $\ppd(4,8\ell-2)\cup\ppd(4,8\ell-6)\cup\ppd(4,4\ell-1)$, we conclude that $|K|$ is divisible by each prime in $\ppd(4,8\ell-2)\cup\ppd(4,8\ell-6)\cup\ppd(4,4\ell-1)$. Then by~\cite[Theorem~A]{LPS1990}, the factorization $B=(A\cap B)K$ implies that $K\cap L=\Sp_{8\ell-2}(4)$.

\underline{Case~1.} Let $\Soc(A)=\Omega_{4\ell}^+(16)$ lie in the second column of row~1 of Table~\ref{tab5}. In this case, $H\cap\Soc(A)=U\times\Sp_{2\ell}(16)$, where $U$ is an nonsolvable subgroup of $\Sp_2(16)$. It follows that $U=\Sp_2(4)$ or $\Sp_2(16)$, and then $H\cap\Soc(A)=\Sp_2(4^c)\times\Sp_{2\ell}(16)$ with $c=1$ or $2$. Since
$$
(H\cap L)/(H\cap\Soc(A))\cong(H\cap L)\Soc(A)/\Soc(A)\leqslant(A\cap L)/\Soc(A)=2^2,
$$
we have $\Sp_2(4^c)\times\Sp_{2\ell}(4)\leqslant H\cap L\leqslant(\Sp_2(4^c)\times\Sp_{2\ell}(4)).2^2$. Thus, $(L,H\cap L,K\cap L)$ lies in row~3 of Table~\ref{tab5}.

\underline{Case~2.} Let $\Soc(A)=\Omega_{4\ell}^+(16)$ lie in the second column of row~4 of Table~\ref{tab5}. In this case, $\ell=3$ and $H\cap\Soc(A)=\Sp_2(4^c)\times\G_2(16)$ with $c=1$ or $2$. As $H\cap\Soc(A)\leqslant H\cap L\leqslant(H\cap\Soc(A)).2^2$, we have $\Sp_2(4^c)\times\G_2(16)\leqslant H\cap L\leqslant(\Sp_2(4^c)\times\G_2(16)).2^2$. Hence $(L,H\cap L,K\cap L)$ lies in row~6 of Table~\ref{tab5}.
\end{proof}

\begin{lemma}\label{lem15}
Suppose $(L,A\cap L,B\cap L)$ lies in \emph{row~6} of \emph{Table~\ref{tab3}}. Then $(L,H\cap L,K\cap L)$ lies in \emph{Table~\ref{tab5}}.
\end{lemma}

\begin{proof}
Here $L=\Omega_8^+(2^f)$ with $f\geqslant2$ and $A\cap L\cong B\cap L=\Sp_6(2^f)$. Since $A=H(A\cap B)$ and $A\cap B\cap L=\G_2(2^f)$ (see~\cite[5.1.15]{LPS1990}), we derive from Hypothesis~\ref{hyp2} that $A\cap L$ is described in the second column of row~7 of Table~\ref{tab1}. It follows that $H\cap L=H\cap(A\cap L)=U\times\Sp_4(2^f)$, where $U$ is an nonsolvable subgroup of $\Sp_2(2^f)$. In particular, $|H\cap L|$ is not divisible by any prime in $\ppd(2,6f)\cup\ppd(2,3f)$. In view of the factorization $B=(A\cap B)K$ we deduce from~\cite[Theorem~A]{LPS1990} that either $K\cap L=\Sp_6(2^f)$, or $K\cap L$ is contained in $\GO_6^+(2^f)$, $\GO_6^-(2^f)$, $\Pa_1[\Sp_6(2^f)]$ or $\N_1[\Sp_6(2^f)]$. For the latter, there exists $r\in\ppd(2,6f)\cup\ppd(2,3f)$ dividing $|L|/|K\cap L|$, which indicates that $r$ divides $|H\cap L|$ by Lemma~\ref{lem8}, a contradiction. Therefore, $K\cap L=\Sp_6(2^f)$, and so $(L,H\cap L,K\cap L)$ lies in row~1 of Table~\ref{tab5}.
\end{proof}

We are now able to complete the proof of Theorem~\ref{thm1}.
\vskip0.1in
\noindent\textit{Proof of Theorem}~\ref{thm1}. Suppose that $G$ is an almost simple group with socle $L$ and $G=HK$, where $H$ has at least two nonsolvable composition factors and $K$ is core-free. We prove Theorem~\ref{thm1} by induction on $|G|$. By Lemma~\ref{lem2} we may assume that $L$ is a classical group. Moreover, by virtue of Lemma~\ref{Embedding} we may let $A$ and $B$ be maximal subgroups of $G$ containing $H$ and $K$, respectively. Thus we have both Hypothesis~\ref{hyp1} and Hypothesis~\ref{hyp2}. Then according to Lemma~\ref{prop2}, either $A$ has at least two nonsolvable composition factors, or the triple $(L,A\cap L,B\cap L)$ lies in Table~\ref{tab3}. If $A$ has at least two nonsolvable composition factors, then Lemma~\ref{prop1} asserts that $(L,H\cap L,K\cap L)$ lies in Table~\ref{tab1} or Table~\ref{tab5}. If $(L,A\cap L,B\cap L)$ lies in Table~\ref{tab3}, then it follows from Lemmas~\ref{lem3}--\ref{lem15} that $(L,H\cap L,K\cap L)$ lies in Table~\ref{tab1} or Table~\ref{tab5}. Hence Theorem~\ref{thm1} is true.
\qed

\end{document}